\documentclass[mnsc]{informs3}

\usepackage{natbib}
 \bibpunct[, ]{(}{)}{,}{a}{}{,}%
 \def\bibfont{\small}%
 \def\newblock{\ }%
\usepackage[colorlinks=true,breaklinks=true,bookmarks=true,urlcolor=blue,
     citecolor=blue,linkcolor=blue,bookmarksopen=false,draft=false]{hyperref}

\def\EMAIL#1{\href{mailto:#1}{#1}}% When hyperref is used, otherwise outcomment 
\TheoremsNumberedThrough     % Preferred (Theorem 1, Lemma 1, Theorem 2)
\EquationsNumberedThrough    % Default: (1), (2), ...
\usepackage{thmtools,thm-restate}
\usepackage{url}
\usepackage[verbose]{wrapfig}
\usepackage{blindtext}
\usepackage{rotating}
\usepackage{enumerate}
\usepackage{amsmath}
\usepackage{float}
\usepackage{verbatim}
\usepackage{amssymb}
\usepackage{color}
\usepackage{bbm}
\usepackage[TABBOTCAP]{subfigure}
\usepackage{tikz}
\usetikzlibrary{shapes}
\usetikzlibrary{arrows}
\usetikzlibrary{topaths}
\usepackage{stmaryrd}
\usepackage{pdfsync}
\usepackage{graphicx}
\newcommand{\set}[1]{\left\{#1\right\}}                     % Set operator
\newcommand{\abs}[1]{\left|#1\right|}                       % Absolute value
\newcommand{\bra}[1]{\left(#1\right)}     

\newcommand{\sidx}[1]{\left\llbracket     #1 \right\rrbracket}  

\newcommand{\Real}{\mathbb R}
\newcommand{\e}{\mathbf{e}}

\newcommand{\Int}{\mathbb Z}

\DeclareMathOperator{\size}{size}
\DeclareMathOperator{\mmc}{mc}

\DeclareMathOperator{\conv}{conv}
\DeclareMathOperator{\ext}{ext}
\DeclareMathOperator{\ray}{ray}

\DeclareMathOperator{\cone}{cone}
\DeclareMathOperator{\aff}{aff}
\DeclareMathOperator{\spann}{span}

\DeclareMathOperator{\gr}{gr}

\usepackage[symbol*]{footmisc}
\DefineFNsymbolsTM*{mylamport*}{%
  \textdagger    \dagger
  \textdaggerdbl \ddagger
  \textsection   \mathsection
  \textparagraph \mathparagraph
  \textbardbl    \|%
  {\textdagger\textdagger}{\dagger\dagger}%
  {\textdaggerdbl\textdaggerdbl}{\ddagger\ddagger}%
  {\textsection\textsection}{\mathsection\mathsection}%
  {\textparagraph\textparagraph}{\mathparagraph\mathparagraph}%
  {\textdagger\textdagger\textdagger}{\dagger\dagger\dagger}%
  {\textdaggerdbl\textdaggerdbl\textdaggerdbl}{\ddagger\ddagger\ddagger}%
  {\textsection\textsection\textsection}%%
    {\mathsection\mathsection\mathsection}%
  {\textparagraph\textparagraph\textparagraph}%%
    {\mathparagraph\mathparagraph\mathparagraph}%
}
\setfnsymbol{mylamport*}
\begin{document}

\TITLE{Embedding Formulations and Complexity for Unions of Polyhedra}
\ARTICLEAUTHORS{%
\AUTHOR{Juan Pablo Vielma}
\AFF{Sloan School of Management, Massachusetts Institute of Technology, \EMAIL{jvielma@mit.edu}}
% Enter all authors
} % end of the block

\ABSTRACT{%

It is well known that selecting a good Mixed Integer Programming (MIP) formulation is crucial for an effective solution with state-of-the art solvers. While best practices and guidelines for constructing good formulations abound, there is rarely a systematic construction leading to the best possible formulation. We introduce embedding formulations and complexity as a new MIP formulation paradigm for systematically constructing formulations for disjunctive constraints that are optimal with respect to size. More specifically, they yield the smallest possible ideal formulation (i.e. one whose LP relaxation has integral extreme points) among all formulations that only use 0-1 auxiliary variables. We use the paradigm to characterize optimal formulations for SOS2 constraints and certain piecewise linear functions of two variables. We also show that the resulting formulations can provide a significant computational advantage over all known formulations for piecewise linear functions.

}%

\MSCCLASS{}
\ORMSCLASS{Primary: ; secondary: }
\HISTORY{29-Apr-2016, 25-Jan-2017, 23-Apr-2017, 18-May-2017}

\maketitle

\section{Introduction}

In its more than 50 years of history  Mixed Integer Programming (MIP) has become an indispensable tool in Operations Research and Management Science. Enormous strides have been made  in the theoretical and  computational issues arising in solving  MIP problems and commercial  MIP solvers  can  solve a wide range of problems  \citep{50book}. One of the reasons for the success of  MIP is its modeling flexibility. For instance, 0-1 MIP can be used to model disjunctive constraints (i.e. the selection over a finite number of alternatives) appearing in a wide range of applications in transportation \citep{croxton2003comparison,roberti2014fixed}, telecommunication \citep{d2013gub} and scheduling \citep{manne1960job,pinedo2012scheduling}. Formulating problems using MIP is often straightforward. However, as most textbooks warn, some care should be taken in constructing MIP formulations as some formulation attributes can severely affect the effectiveness of solvers. Good formulations can be obtained following simple guidelines, but more elaborate techniques can provide a significant computational advantage \citep{Mixed-Integer-Linear-Programming-Formulation-Techniques}. For instance, consider the classical \emph{Special Ordered Sets of Type 2 (SOS2)} introduced by  \cite{Beale70}. SOS2 constraints on variables $\lambda\in \mathbb{R}^{n+1}$ require that (i) $\sum\nolimits_{i=1}^{n+1} \lambda_i=1$ and $\lambda_i\geq 0$ for all $i\in \set{1,\ldots,n+1}$, and (ii) at most two $\lambda_i$ variables can be non-zero at the same time and
if $\lambda_i>0$ and $\lambda_j>0$, then they must be adjacent variables (i.e. $\abs{i-j}\leq 1$). A textbook formulation for SOS2 constraints for $n=4$ is given by 
\begin{subequations}\label{CC1d}
\begin{alignat}{3}
\sum\nolimits_{j=1}^5 \lambda_j=1,\quad \lambda_j \geq 0 \quad \forall j\in \set{1,\ldots,5}, \quad\sum\nolimits_{i=1}^4 y_i=1,\quad y\in \set{0,1}^4\\
\lambda_1 \leq y_1,\quad \lambda_2\leq y_1+y_2,\quad \lambda_3 \leq y_2+y_3,\quad \lambda_4\leq y_3+y_4,\quad \lambda_5\leq y_4.\label{CC1dg}
\end{alignat}
\end{subequations}
To evaluate this formulation we consider two formulation attributes: formulation size and strength. For size we use the standard measure corresponding to the smallest number of linear inequalities needed to describe the formulation (where we count an equation as two inequalities). Formulation strength can be evaluated in many ways  \cite[Section 2.2]{Mixed-Integer-Linear-Programming-Formulation-Techniques};  we choose to check if the formulation satisfies the following property.
\begin{definition}[Ideal Formulation] Let   $Ax+By \leq b,\quad y\in \mathbb{Z}^k$ be a MIP formulation\footnote{Where $A$, $B$ and $b$ are appropriately sized rational matrices and vector.}. The Linear Programming (LP) relaxation of this formulation is the polyhedron described by $Ax+By \leq b$, which for simplicity we assume has at least one extreme point or basic feasible solution\footnote{For a case in which this assumption is relaxed see \cite{vielma2017small}.}. We say the MIP formulation is \emph{ideal} if and only if all the extreme points of its LP relaxation satisfy the integrality constraints $y\in \mathbb{Z}^k$.
\end{definition}

Ideal formulations are in a sense \emph{strongest possible} \cite[Section 2.2]{Mixed-Integer-Linear-Programming-Formulation-Techniques} and usually outperform similarly sized non-ideal formulations. Formulation \eqref{CC1d} fares well with regards to size, as its version for general $n$ only requires $2n+6$ inequalities, with $4$ of these coming from equations (the bounds on the $y$ variables are not actually needed). However, this small size comes at the cost of the formulation not being ideal. Indeed, for \eqref{CC1d} to be ideal every extreme point $\bra{\lambda,y}$ of its LP relaxation must satisfy $y\in \mathbb{Z}^4$. The LP relaxation of \eqref{CC1d} is obtained by replaxing $y\in \set{0,1}^4$ by $0\leq y_i\leq 1$ for all $i\in \set{1,\ldots,4}$ and we can check that $\lambda=(1/2,1/2,0,0,0)$ and $y=(1/2,0,1/2,0)$ is an extreme point of this resulting set. Fortunately, \cite{padberg00} showed how to strengthen \eqref{CC1d} to an ideal formulation without increasing the total number of inequalities. For $n=4$ the resulting formulation is given by 
\begin{subequations}\label{CC1dstrong}
\begin{alignat}{3}
\sum\nolimits_{j=1}^5 \lambda_j=1,\quad \lambda_1\geq 0,\quad\lambda_5 \geq 0, \quad\sum\nolimits_{i=1}^4 y_i=1,\quad y\in \set{0,1}^4\\
\lambda_1\leq y_1\leq \lambda_1+\lambda_2\leq y_1+y_2\leq \lambda_1+\lambda_2+\lambda_3\leq y_1+y_2+y_3\leq \lambda_1+\lambda_2+\lambda_3+\lambda_4.\label{CC1dstrongg}
\end{alignat}
\end{subequations}
While both formulations have nearly the same size, the majority of the inequalities of \eqref{CC1dstrong} correspond to general inequalities \eqref{CC1dstrongg}, which are not just variable bounds. In contrast, only half the inequalities of \eqref{CC1d} are general inequalities (see Table~\ref{formsizetable} and Corollary~\ref{unaryformulationsos2}). General inequalities usually have a stronger computational impact than variable bounds as the later can be treated implicitly by LP and MIP solvers. Hence the larger number of general inequalities of  \eqref{CC1dstrong} can cancel its advantage from being ideal and indeed it is often  computationally outperformed by non-ideal formulation \eqref{CC1d}. A solution to this issue can be found in an advanced  formulation technique introduced by  \cite{Modeling-Disjunctive-Constraints-FULL}. For $n=4$ this technique yields the formulation for SOS2 given by
\begin{subequations}\label{LogCC1dstrong}
\begin{alignat}{3}
\sum\nolimits_{j=1}^5 \lambda_j=1,\quad \lambda_j \geq 0 \quad \forall j\in \set{1,\ldots,5}, \quad y\in \set{0,1}^2\\
\lambda_1+\lambda_5\leq  1-y_1,\quad \lambda_3\leq y_1,\quad \lambda_4+\lambda_5\leq 1-y_2,\quad \lambda_1+\lambda_2\leq y_2.
\end{alignat}
\end{subequations}
This formulation is also ideal, but its version for general $n$ requires around $2\lceil \log_2 n\rceil+n+2$ inequalities and only $2\lceil \log_2 n\rceil+2$ of these are general inequalities (see Corollary~\ref{sosub}). This allows formulation \eqref{LogCC1dstrong} to have a significant computational advantage over \eqref{CC1d}, \eqref{CC1dstrong} and all known formulations for SOS2 \citep{Mixed-Integer-Models-for-Nonseparable,Modeling-Disjunctive-Constraints-FULL}.  An even more dramatic issue arises if we consider a 2-dimensional generalization of SOS2 constraints used to model piecewise linear functions of two variables \citep{lee01,Mixed-Integer-Models-for-Nonseparable}. As detailed in Table~\ref{formsizetable}, the cost of going from the non-ideal 2-dimensional generalization of  \eqref{CC1d} to the ideal 2-dimensional generalization of \eqref{CC1dstrong} is a significant increase in the number of inequalities, particularly general inequalities. However, the 2-dimensional generalization of formulation \eqref{LogCC1dstrong} still has a linear number of inequalities and a logarithmic number of general inequalities. This again gives it a significant computational advantage over all known formulations for piecewise linear functions of two variables \citep{Mixed-Integer-Models-for-Nonseparable,Modeling-Disjunctive-Constraints-FULL}.
\begin{table}[htpb]
\begin{center}
\bgroup
\def\arraystretch{2}%
\begin{tabular}{c|ll|ll}
& \multicolumn{2} {l|}{Traditional SOS2}& \multicolumn{2} {l}{2-D Generalization}\\
Formulation & General Inequalities &Bounds& General Inequalities &Bounds\\
\hline
\eqref{CC1d}& $n+1$ &$n+1$ & $\bra{\sqrt{n/2}+1}^2$ & $\bra{\sqrt{n/2}+1}^2$\\[2ex]
\eqref{CC1dstrong}& $2(n-1)$ &$2$& $\dbinom{2\sqrt{n/2}}{\sqrt{n/2}}$ & $\bra{\sqrt{n/2}+1}^2$\\[2ex]
\eqref{LogCC1dstrong} & $2\lceil \log_2 n\rceil$ &$[n-1,n+1]$& $4 \left\lceil \log_2 \sqrt{n/2}\right\rceil+ 2$ & $\bra{\sqrt{n/2}+1}^2$\\
\end{tabular}
\egroup
  \end{center}
\caption{Sizes of Formulations for SOS2 Constraints and its 2-dimensional Generalization. We omit the number of inequalities from equations, which is $4$ for \eqref{CC1d} and \eqref{CC1dstrong}, and $2$ for \eqref{LogCC1dstrong}.}\label{formsizetable}
\end{table}
 
The only dissadvantage of formulation \eqref{LogCC1dstrong} is its increase in complexity. In particular, while for formulations \eqref{CC1d} and  \eqref{CC1dstrong} we can easily interpret the role of the $0$-$1$ variables ($y_i=1$ if and only if $\lambda_i$ and $\lambda_{i+1}$ can be non-zero at the same time), the role of the $0$-$1$ variables is not so clear for \eqref{LogCC1dstrong} (for instance it uses two $0$-$1$ variables instead of four). This increase in complexity makes it hard to generalize \eqref{LogCC1dstrong} to other constraints. In fact, the 2-dimensional generalization of  \eqref{LogCC1dstrong} only works for very specific piecewise linear functions, while the 2-dimensional generalizations of \eqref{CC1d} and \eqref{CC1dstrong} work for a wide range of piecewise linear functions.

In this paper we propose a new MIP formulation paradigm that should allow extending the success of formulation \eqref{LogCC1dstrong} to a wide range of applications. In particular, this paradigm can construct ideal formulations for any disjunctive constraint that requires a set of variables to be in the union of a finite number of polyhedra with mild technical requirements. We denote the formulations obtained through this paradigm \emph{embedding formulations} as they are based on a geometric construction that \emph{embeds} the disjunctive constraints into a higher dimensional space that contains both the original variables in the constraint (e.g. the $\lambda$ variables for SOS2 constraints) and the $0$-$1$ variables of the formulation (e.g. the $y$ variables). One characteristic of embedding formulations is allowing a flexible use of $0$-$1$ variables that include both traditional uses such as in \eqref{CC1d} and  \eqref{CC1dstrong},  and more complex uses such as in \eqref{LogCC1dstrong}. This flexibility is the key to replicating the success of \eqref{LogCC1dstrong} as we show that the size of an embedding formulations can be extremely sensitive to the specific use of $0$-$1$ variables. For this reason we also study the  size of the smallest embedding formulation for a disjunctive constraint when we consider all possible uses of $0$-$1$ variables. We denote this the embedding complexity of the associated union of polyhedra. 
This complexity measure has theoretical interest on its own, but can also be used to evaluate the potential for improvement of existing formulations. For instance, studying this complexity allows us to show that \eqref{LogCC1dstrong} is (nearly) optimal with regard to size. To the best of our knowledge, this result is the first lower bound on sizes of \emph{mixed} integer formulations and one of the few results for any class of integer programming formulations. The only other similar results we are aware of are the techniques to show lower bounds on sizes of combinatorial and \emph{pure} integer formulations introduced by \cite{KaibelW15b,weltge2015sizes}. Finally, we show how the embedding formulation paradigm can be used to generalize the 2-dimensional version of \eqref{LogCC1dstrong} to a wider range of piecewise linear functions than what was considered in  \cite{Modeling-Disjunctive-Constraints-FULL}. We also show how the resulting formulations can significantly outperform all other known formulations for piecewise linear functions of two variables. This generalization is based on the computational calculation of a convex hull associated to the embedding formulation. To the best of our knowledge, this is the first example of a computational construction of an effective MIP formulation. 

Throughout  the paper we use the following notation. For a set $S\subseteq \Real^d$ we let $\conv\bra{S}$, $\aff\bra{S}$, $\spann\bra{S}$ and $\dim\bra{S}$ be the convex hull, affine hull, linear span and the dimension of $S$ respectively. For a polyhedron $P\subseteq \Real^d$ we let $\ext\bra{P}$ and $\ray\bra{P}$ be the set of extreme points and extreme rays of $P$. We also let $P_\infty$ be the recession cone of $P$. Given two vectors $a, b\in \Real^V$ for a finite index set $V$ we let $a\cdot b=\sum\nolimits_{v\in V} a_v b_v$ be the inner product between $a$ and $b$. We also let ${\bf 0}\in \Real^V$ be the vector of all zeros and  $\e^v\in \Real^V$ be the unit vector such that $\e^v_u=1$ if $u=v$ and $\e^v_u=0$ otherwise (The index set $V$ is often evident form the context so we omit it in this notation). Finally, we let $\sidx{n}:=\set{1,\ldots,n}$, $\sidx{a,b}:=\set{a,a+1,\ldots, b-1,b}$ and $\mathbb{Q}$ be the set of rational numbers.

\section{Geometric Construction of Formulations}\label{formulationsec}

\begin{figure}[htpb]
  \begin{center}
  \includegraphics[scale=0.35]{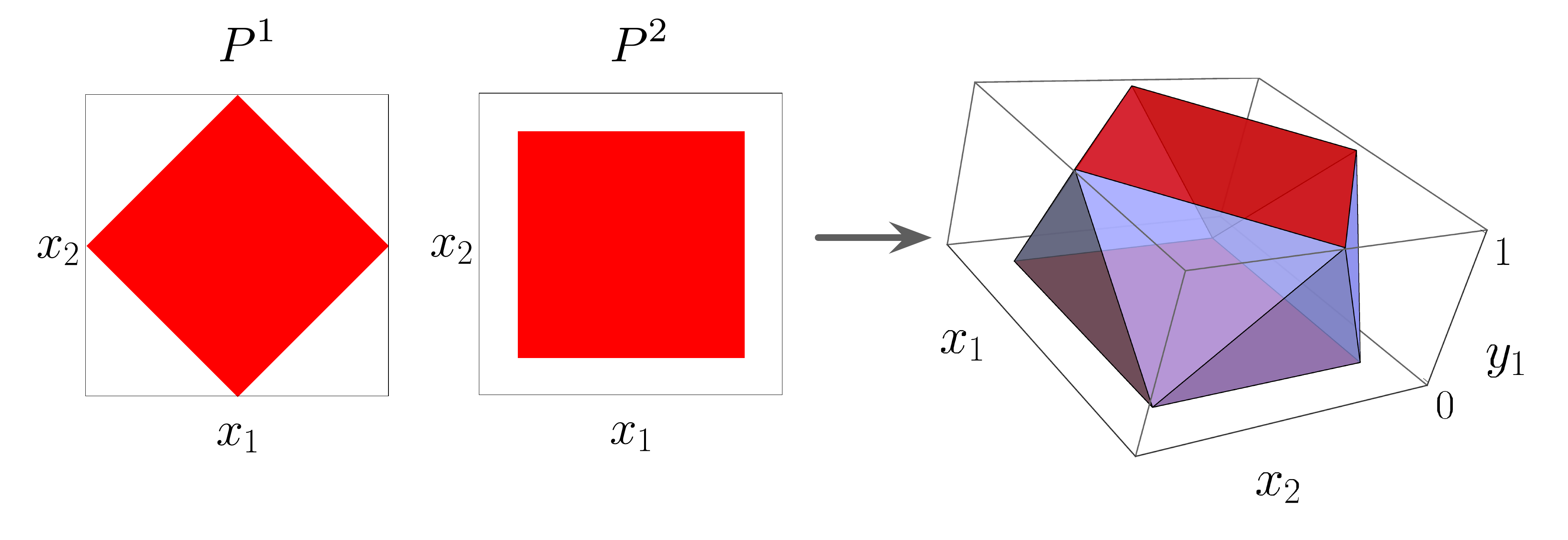}
  \end{center}
  \caption{Two polyhedra (left in red), their embedding into a space with one $0$-$1$ variable (right in red) and the convex hull of this embedding (right in light blue).}\label{embeddingillust}
  \end{figure}

We consider MIP formulations for the disjunctive constraint \begin{equation}\label{disjunction}
x\in \bigcup\nolimits_{i=1}^n P^i
\end{equation}
where $\mathcal{P}:=\bra{P^i}_{i=1}^n$ is a finite family of polyhedra in $\Real^d$ that satisfy the following assumption.
\begin{assumption}\label{ass1}The family of polyhedra $\mathcal{P}:=\bra{P^i}_{i=1}^n$ is such that  $P^i$ is a rational polyhedron and $\ext\bra{P^i}\neq \emptyset$ for all $i\in \sidx{n}$, and 
 $P^i_\infty=P^j_\infty$ for all $i,j\in \sidx{n}$.
\end{assumption}
Disjunctive constraint \eqref{disjunction} is exactly the type of practical\footnote{We can also consider unions of polyhedra without extreme points, but this is uncommon in practice.} constraints that can be modeled using $0$-$1$ MIP \citep{springerlink:10.1007/BFb0121015} and MIP formulations for it can be constructed using a wide range of techniques \citep{Mixed-Integer-Linear-Programming-Formulation-Techniques}.
For instance, if the polyhedra are described by linear inequalities, the following result from \cite{springerlink:10.1007/BFb0121015,balas85} gives an ideal and small formulation for \eqref{disjunction}. 
\begin{theorem}\label{BTheo}  Let $\mathcal{P}:=\bra{P^i}_{i=1}^n$ be a finite family of polyhedra in $\Real^d$ satisfying Assumption~\ref{ass1}. Furthermore, for each $i\in\sidx{n}$, let $A^i\in \mathbb{Q}^{m_i\times d}$ and $b^i\in \mathbb{Q}^{m_i}$ be such that $P^i=\set{x\in \Real^d\,:\, A^i x\leq b^i}$. Then an ideal formulation of \eqref{disjunction} is given by 
\begin{alignat}{3}\label{BalasForm}
A^i x^i\leq b^i y_i\quad\forall i\in \sidx{n},\quad x=\sum\nolimits_{i=1}^n x^i,\quad \sum\nolimits_{i=1}^n y_i=1,\quad y\in \set{0,1}^n.
\end{alignat}
\end{theorem}
Formulation \eqref{BalasForm} achieves the feat of being ideal and small through the use of continuous auxiliary variables $x^i\in \Real^d$, which copy the original variables $x$ for each polyhedron. However, these continuous auxiliary variables can take away the potential computational advantage of \eqref{BalasForm}, particularly when small and ideal formulations without the continuous auxiliary variables are available (e.g. formulation \eqref{LogCC1dstrong} for SOS2 constraints). Unfortunately, there is no known general  formulation for disjunctive constraints that is small, ideal and  does not use these continuous auxiliary variables. To remedy this, we  present a general geometric procedure to construct formulations without continuous auxiliary variables that is always ideal and can sometimes yield small formulations. The basic steps of this procedure are depicted  in Figure~\ref{embeddingillust} for the polyhedra $P^1$ and $P^2$ in the left of the figure. The first step embeds the polyhedra into a space that contains an additional $0$-$1$ variable $y_1$ by converting the disjunction from $S:=P^1\cup P^2$ to  $S^+:=\bra{P^1\times \set{1}} \cup \bra{P^2\times \set{0}}$. The second step takes the convex hull of this embedding to obtain $Q=\conv\bra{S^+}$. By construction we have that $S^+=Q\cap \bra{\mathbb{R}^2\times \mathbb{Z}}$ and hence the projection of $Q\cap \bra{\mathbb{R}^2\times \mathbb{Z}}$ onto the $x$ variables is equal to $S$. Furthermore, by construction we also have that the extreme points of $Q$ have an integral $y_1$ component. Hence, an ideal formulation of $S$ is given by $\bra{x,y_1}\in Q=\conv\bra{S^+}$ and $y_1\in \mathbb{Z}$. Finally, because $Q$ is a rational polyhedron, this formulation contains only linear inequalities with rational coefficients and the integrality constraint on $y_1$. 

\subsection{Embedding Formulations}

One possible generalization of the embedding procedure in Figure~\ref{embeddingillust} to more than two polyhedra is to embed $\bigcup_{i=1}^n P^i$ into $\bigcup_{i=1}^n \bra{ P^i \times \set{\e^i}}$ where $\e^i$ is the $i$-th unit vector. This procedure is known as the \emph{Cayley trick} or \emph{Cayley Embedding} and is used to study Minkowski sums of polyhedra (e.g. \cite{caytrick,karavelas2013maximum,WeibelPhd}). We here consider a further generalization that pairs the polyhedra $P^i$ with any set of pairwise disjoint $0$-$1$ vectors instead of just the unit vectors $\e^i$. The following proposition and its proof shows that the convex hull of the resulting embedding is a rational polyhedron that can be used to construct an ideal formulation of $x\in \bigcup_{i=1}^n P^i$. 

\begin{proposition}\label{firstprop} Let $\mathcal{P}:=\bra{P^i}_{i=1}^n$ be a finite family of polyhedra in $\Real^d$ satisfying Assumption~\ref{ass1}, $k\geq \left\lceil \log_2 n\right\rceil$,  
$H:=\bra{h^i}_{i=1}^n\in \mathcal{H}_k(n):=\set{\bra{h^i}_{i=1}^n\,:\, h^i\in \set{0,1}^k \quad \forall i\in \sidx{n},\quad h^i\neq h^j \quad \forall i\neq j}$
be a family of pairwise distinct $0$-$1$ vectors indexed in the same way as the polyhedra and
\begin{equation}\label{embconv}
Q\bra{\mathcal{P},H}:=\conv\bra{\bigcup\nolimits_{i=1}^n P^i\times \set{h^i} }.
\end{equation}
 Then 
\begin{enumerate}[(i)]
\item\label{rationalLPprop} $Q\bra{\mathcal{P},H}$ is a rational polyhedron,
\item\label{validprop} $ \bra{x,y}\in Q\bra{\mathcal{P},H} \cap \bra{\Real^d\times\Int^k} \quad \Leftrightarrow \quad\exists i\in \sidx{n} \text{ s.t. } y=h^i\;\wedge\; x\in P^i,$ and
\item\label{idealprop} $\ext\bra{Q\bra{\mathcal{P},H}}\subseteq \Real^d\times \set{0,1}^k$.
\end{enumerate}
\end{proposition}
\proof{Proof.}
By Assumption~\ref{ass1},  $\set{P^i\times \set{h^i}}_{i=1}^n$ is a finite family of non-empty polyhedra with identical recession cones. Then, by Lemma 4.41 and Corollary 4.44  in \cite{conforti2014integer}, $Q\bra{\mathcal{P},H}$ is a rational polyhedron and
$Q\bra{\mathcal{P},H}=\conv\bra{\bigcup_{i=1}^n \bigcup_{v\in \ext\bra{P^i}  }\set{v\times h^i}}+\cone\bra{ \bigcup_{r\in \ray\bra{P^1}  }\set{r\times {\bf 0}}}$.
This shows the first two properties of $Q\bra{\mathcal{P},H}$. The last follows by additionally noting that the $h^i$'s are distinct extreme points of $[0,1]^k$.
\Halmos\endproof

Similarly to  the example depicted in Figure~\ref{embeddingillust} we obtain a formulation from $Q\bra{\mathcal{P},H}$ by interpreting $H$ as possible values of a $k$-dimensional  $0$-$1$ variable $y$.
\begin{corollary}[Embedding Formulation]\label{embeddingformdef}
 Let $\mathcal{P}:=\bra{P^i}_{i=1}^n$ be a finite family of polyhedra in $\Real^d$ satisfying Assumption~\ref{ass1}, $k\geq \left\lceil \log_2 n\right\rceil$ and $H\in \mathcal{H}_k(n)$. Then an ideal formulation of $x\in \bigcup_{i=1}^n P^i$ is given by 
 \begin{equation}\label{embformeq}
  \bra{x,y}\in Q\bra{\mathcal{P},H},\quad y\in \mathbb{Z}^k.
  \end{equation}
 Furthermore, besides the integrality constraints on $y$, this formulation only includes linear inequalities with rational coefficients. We refer to \eqref{embformeq} as the \emph{embedding formulation} of $\mathcal{P}$ associated to $H$, to $H$ as the \emph{encoding} of the formulation and to $Q\bra{\mathcal{P},H}$ as the \emph{LP relaxation} of the formulation.
\end{corollary}
\proof{Proof.} The formulation is valid because of property \eqref{validprop} in Proposition~\ref{firstprop}, includes only linear inequalities with rational coefficients because of property \eqref{rationalLPprop} in Proposition~\ref{firstprop}, and is ideal because of property \eqref{idealprop} in Proposition~\ref{firstprop}.  \Halmos\endproof

Encoding $H$ can have a strong impact in the size of $Q\bra{\mathcal{P},H}$, but different encodings can yield formulations of the same size or even equivalent polyhedra $Q\bra{\mathcal{P},H}$. The following straightforward lemma shows one such possible equivalency. 
\begin{lemma}\label{isolemma}
Let $\mathcal{P}:=\bra{P^i}_{i=1}^n$ be a finite family of polyhedra in $\Real^d$ satisfying Assumption~\ref{ass1}, $k_1,k_2\geq \left\lceil \log_2 n\right\rceil$, $H\in \mathcal{H}_{k_1}(n)$ and $G\in \mathcal{H}_{k_2}(n)$. If there exists an affine map  $A:\Real^{k_1}\to \Real^{k_2}$ such that $A$ is a bijection between $\conv\bra{H}$ and $\conv\bra{G}$ then $Q\bra{\mathcal{P},H}$ is affinely isomorphic  to $Q\bra{\mathcal{P},G}$.

In particular, if $H:=\bra{h^i}_{i=1}^n$  where $h^i=\e^i$ for all $i\in \sidx{n}$, and  $G:=\set{g^i}_{i=1}^n$  where $g^i=\e^{\pi(i)}$ for all $i\in \sidx{n}$ and $\pi:\sidx{n}\to \sidx{n}$ is a permutation, then 
\[\bra{x,y}\in Q\bra{\mathcal{P},G} \quad \Leftrightarrow \quad \bra{x,y^\pi}\in Q\bra{\mathcal{P},H}  \]
where $y^\pi_{\pi\bra{i}}=y_i$ for all $i\in \sidx{n}$.
\end{lemma}

Lemma~\ref{isolemma} shows that when the encoding uses $n$ unit vectors $\e^i\in \set{0,1}^n$ the specific order or polytope-vector assignment in the embedding is inconsequential. Hence, when an encoding of this form is used we assume the unit vectors are assigned in their natural order and refer to the resulting embedding formulation as the \emph{unary encoded} formulation as this encoding can be interpreted as a \emph{unary encoding} of the selection among the polyhedra. A completely different class of encodings are obtained when  $n$ is a power of $2$ and  $\bra{h^i}_{i=1}^n=\set{0,1}^k$ for $k=\log_2 n$. This case can be interpreted as a \emph{binary encoding} of the selection among the polyhedra and corresponds to the encodings with the smallest number of components or \emph{bits}. For this reason we refer to embedding formulations resulting from such encodings as \emph{binary encoded} formulations. Unlike unary encoded formulations, in Section~\ref{SOS2} we show that permuting the order of a binary encoding can lead to  binary encoded formulations of significantly different sizes. This potential size variability over binary and other encodings motivates the following complexity measure for unions of polyhedra, which quantifies the size of its smallest embedding formulation.

\begin{definition}[Embedding Complexity] For a polyhedron $Q$ let $\size(Q)$ be equal to the minimum number of inequalities needed to describe $Q$ (equations are counted as two inequalities).  Then, for a family of polyhedra $\mathcal{P}:=\bra{P^i}_{i=1}^n$ satisfying Assumption~\ref{ass1} we let its  \emph{embedding complexity} be 
\[\mmc\bra{\mathcal{P}}:=\min\set{\size\bra{Q\bra{\mathcal{P},H}}\,:\,  H\in \bigcup\nolimits_{k\geq \left\lceil \log_2 n\right\rceil}  \mathcal{H}_k(n) }.\]
\end{definition}

Constructing an embedding formulation even for a fixed $H$ requires a potentially costly convex hull calculation. Calculating the embedding complexity has the added difficulty of minimizing the size of these convex hulls over all possible encodings. Fortunately, as we show in the following sections it is sometimes possible to give tight bounds on the embedding complexity of specially structured disjunctions such as SOS2 constraints. We also show how we can computationally construct embedding formulations for piecewise linear functions that can provide a significant computational advantage. In both these cases the unions of polyhedra considered have a special combinatorial structure (specifically, all polyhedra are faces of a fixed simplex). For an example of how embedding formulations can lead to a computational advantage for disjunctions without this special structure we refer the reader to \cite{huchette2016strong}.

\section{Bounds on embedding complexity for SOS2 constraints}\label{SOS2}

SOS2 constraints can be posed as a disjunctive constraint of the form \eqref{disjunction} as follows.
\begin{definition}
Let $\Delta^{n+1}:=\set{\lambda\in \Real^{n+1}_+\,:\, \sum_{i=1}^{n+1} \lambda_j=1}$. A family of polyhedra $\mathcal{P}:=\bra{P^i}_{i=1}^n$ in $\Real^{n+1}$ is the special ordered sets of type 2 or SOS2 constraint on $\Delta^{n+1}$ if and only if
\[ P^i:=\set{\lambda\in \Delta^{n+1}\,:\,  \lambda_j\leq 0\quad \forall j\not\in \set{i,i+1}}=\conv\bra{\set{\e^i,\e^{i+1}}}.\]
\end{definition}

To fully characterize embedding formulations for SOS2 constraints we need to consider the possibility of an encoding not being full dimensional (i.e. $H\in\mathcal{H}_k(n)$ and $\dim(H)<k$), which we handle through the following definition.
\begin{definition}For $H:=\bra{h^i}_{i=1}^n$, let   $L\bra{H}:=\aff\bra{H}-h^1$ be the linear space parallel to the affine hull of $H$. 
\end{definition}

The following proposition shows that the size and structure of the inequalities of $ Q\bra{\mathcal{P},H}$ have a relatively simple description for SOS2 constraints. 
	
\begin{restatable}{proposition}{sosprophyperprop}\label{sos1prophyper}Let $\mathcal{P}:=\bra{P^i}_{i=1}^n$ be the SOS2 constraint on $\Delta^{n+1}$, $H:=\bra{h^i}_{i=1}^n\in \mathcal{H}_k(n)$, $c^i=h^{i+1}-h^i$ for $i\in \sidx{n-1}$ and for $b\in L(H)\setminus \set{\bf 0}$ let  $M(b):=\set{y\in L(H)\,:\, b\cdot y=0}$ be the linear (or central) hyperplane defined by $b$ in $L(H)$.
Finally, let $\set{b^l}_{l=1}^L\subseteq L(H)\setminus \set{\bf 0}$ be  such that $\set{M\bra{b^l}}_{l=1}^{L}$  is the set of linear hyperplanes spanned by $\bra{c^i}_{i=1}^{n-1}$ in $L(H)$ and $J=\set{1,n+1}\cup \set{j\in \sidx{2,n}\,:\, L\bra{H}=\spann\bra{\set{c^i}_{i\in \sidx{n-1}\setminus\set{j-1} }}}  $. Then $\bra{\lambda,y}\in Q\bra{\mathcal{P},H}$ if and only if
\begin{subequations}\label{SOS2char}
\begin{alignat}{3}
\sum\nolimits_{j=1}^{n+1} \lambda_j =1,\quad y&\in \aff\bra{H}\label{sos2equalities}\\
\sum\nolimits_{j=1}^{n+1} \min\set{{b^l\cdot h^j} ,{b^l\cdot h^{j-1}}}\lambda_j \leq b^l\cdot y &\leq \sum\nolimits_{j=1}^{n+1} \max\set{{b^l\cdot h^j} ,{b^l\cdot h^{j-1}}}\lambda_j &\quad&\forall l \in \sidx{L}\label{sos2generalfacets}\\
\lambda_j&\geq 0 &\quad& \forall j\in J\label{sos2boundfacets},
\end{alignat}
\end{subequations} 
where we let $h^0=h^1$ and $h^{n+1}=h^n$. Furthermore, $\dim\bra{Q\bra{\mathcal{P},H}}=n+\dim\bra{H}$, equations \eqref{sos2equalities} precisely describe $\aff\bra{Q\bra{\mathcal{P},H}}$ and all inequalities in \eqref{sos2generalfacets} and \eqref{sos2boundfacets} are facet defining for $Q\bra{\mathcal{P},H}$. Finally, none of the facets defined by \eqref{sos2generalfacets} can be defined by a variable bound on $\lambda$.
\end{restatable}

The proof of Proposition~\ref{sos1prophyper} is slightly technical and follows from a complete characterization of the facial structure of $ Q\bra{\mathcal{P},H}$ for SOS2 constraints, so we postpone it to Section~\ref{proofofsos1prophyper}.

One notable characteristic of \eqref{SOS2char} is the simple form of the inequalities. However, this simplicity hides the possibility of discovering rather complicated formulations whose validity is not straightforward. Indeed, combining Proposition~\ref{sos1prophyper} and Corollary~\ref{embeddingformdef} we have that a valid (and ideal) formulation for SOS2 constraints is $\bra{\lambda,y}$ satisfy \eqref{SOS2char} and $y\in \mathbb{Z}^k$. We can easily check that $\bra{\lambda,y}=\bra{\e^j,h^i}$ is feasible for this formulation if $j\in\set{i,i+1}$. The converse can sometimes be easily verified for encodings with a specific structure, but is not evident from generic formulation \eqref{SOS2char}. To complicate matters further, this converse does not directly imply infeasibility of $\bra{\lambda,y}=\bra{1/2}\bra{\e^j,h^i}+\bra{1/2}\bra{\e^l,h^i}$ for $j\in\set{i,i+1}$ and $l\notin\set{i,i+1}$, which is also required for validity of the formulation. In Section~\ref{examplesection} give an example of how even for fixed $n$ and $H$ it may not be straightforward to verify validity of the formulation based on \eqref{SOS2char} \emph{directly} (i.e. without using Corollary~\ref{embeddingformdef} and the fact that \eqref{SOS2char} describes $ Q\bra{\mathcal{P},H}$). This example illustrates how the embedding formulation procedure can be useful to discover or construct unexpected formulations whose validity is not immediately evident. 

Besides constructing unexpected formulations, Proposition~\ref{sos1prophyper} allows us to recover existing formulations and determine that one of these is the smallest embedding formulation for SOS2 constraints. Characterization  \eqref{SOS2char} of 
$ Q\bra{\mathcal{P},H}$ considers two classes of facets: those defined by variable bounds \eqref{sos2boundfacets} and those defined by more general inequalities \eqref{sos2generalfacets}. As noted in the introduction the computational cost of variable bounds tends to be smaller so we refine our notion of formulation size to consider this in our analysis. 
\begin{definition} Let $\mathcal{P}:=\bra{P^i}_{i=1}^n$ be the SOS2 constraint on $\Delta^{n+1}$ and $H:=\bra{h^i}_{i=1}^n\in \mathcal{H}_k(n)$. We let $\size_G\bra{Q\bra{\mathcal{P},H}}$ and $\size_B\bra{Q\bra{\mathcal{P},H}}$ be the number of inequalities in \eqref{sos2generalfacets} and \eqref{sos2boundfacets} respectively so that $\size\bra{Q\bra{\mathcal{P},H}}=\size_G\bra{Q\bra{\mathcal{P},H}}+\size_B \bra{Q\bra{\mathcal{P},H}}+2\bra{1+k-\dim\bra{H}}$. Finally, we let 
$\mmc_G\bra{\mathcal{P}}:=\min\set{\size_G\bra{Q\bra{\mathcal{P},H}}\,:\,  H\in \bigcup\nolimits_{k\geq \left\lceil \log_2 n\right\rceil}  \mathcal{H}_k(n) }$.
\end{definition}

The first formulation we can recover with Proposition~\ref{sos1prophyper} is the ideal formulation  introduced in  \cite{padberg00} that was illustrated in \eqref{CC1dstrong}.
\begin{corollary}\label{unaryformulationsos2}Let $\mathcal{P}:=\bra{P^i}_{i=1}^n$ be the SOS2 constraint on $\Delta^{n+1}$ and  $H$ be the unary encoding, then  $\bra{\lambda,y}\in Q\bra{\mathcal{P},H}$ if and only if
\begin{subequations}\label{padbergform2}
\begin{alignat}{3}
\sum\nolimits_{j=1}^{n+1} \lambda_j=1, \quad\sum\nolimits_{i=1}^n y_i&=1,\quad y\in \set{0,1}^n,
&\quad& \\
\sum\nolimits_{j=l+2}^{n+1} \lambda_j \leq \sum\nolimits_{i=l+1}^n y_i,\quad  
\sum\nolimits_{j=1}^{l} \lambda_j& \leq \sum\nolimits_{i=1}^l y_i &\quad&\forall l \in \sidx{n-1}\label{unaryineq1}\\
\lambda_1\geq 0,\quad \lambda_{n+1}&\geq 0.\label{unaryineq2}
\end{alignat}
\end{subequations}
Furthermore, all inequalities in \eqref{unaryineq1}--\eqref{unaryineq2} are facet defining and hence $\size\bra{Q\bra{\mathcal{P},H}}=2n+4$, $\size_G\bra{Q\bra{\mathcal{P},H}}=2(n-1)$ and $\size_B\bra{Q\bra{\mathcal{P},H}}=2$.
\end{corollary}
\proof{Proof.}
For the unary encoding we have $c^i=\e^{i+1}-\e^i$ for all $i\in \sidx{n-1}$ and hence $L(H)=\set{y\in \Real^{n}\,:\,\sum_{j=1}^{n}y_j=0}$. Furthermore, the set of hyperplanes spanned by $\bra{c^i}_{i=1}^{n-1}$ is $\set{M\bra{b^l}}_{l=1}^{L}$ for $L=n-1$ and $b^l=\bra{n-l}\sum_{i=1}^l \e^i - l\sum_{i=l+1}^n \e^i$ for each $l\in \sidx{L}$. Finally, the elements of $\bra{c^i}_{i=1}^{n-1}$ are linearly independent and hence $J=\set{1,n+1}$ and  \eqref{SOS2char} for the unary encoding becomes 
\begin{subequations}\label{padbergform}
\begin{alignat}{3}
\sum\nolimits_{j=1}^{n+1} \lambda_j=1, \quad\sum\nolimits_{i=1}^n y_i=1,\quad &y\in \set{0,1}^n\\
\label{padberg1}(n-l){\sum\nolimits_{j=1}^l \lambda_j}-l{\sum\nolimits_{j=l+1}^{n+1}\lambda_j}&\leq (n-l){\sum\nolimits_{i=1}^l y_i}-l{\sum\nolimits_{i=l+1}^n y_i}&\quad& \forall l\in \sidx{n-1}\\
\label{padberg2}(n-l){\sum\nolimits_{i=1}^l y_i}-l{\sum\nolimits_{i=l+1}^n y_i}&\leq (n-l){\sum\nolimits_{j=1}^{l+1} \lambda_j}-l{\sum\nolimits_{j=l+2}^{n+1}\lambda_j}
&\quad& \forall l\in \sidx{n-1}.\\
\lambda_1\geq 0,\quad \lambda_{n+1}&\geq 0.
\end{alignat}
\end{subequations}
If we subtract $n-l$ times the implied equation $\sum_{i=1}^n y_i=\sum_{j=1}^{n+1} \lambda_j$ from \eqref{padberg2} and divide by $n$, and add $l$ times this same implied equation to \eqref{padberg1} and divide by $n$ we have that \eqref{padbergform} is equivalent to \eqref{padbergform2}.
\Halmos\endproof

We can also use  Proposition~\ref{sos1prophyper} to recover and generalize the logarithmic formulation from \cite{Modeling-Disjunctive-Constraints-FULL} that we illustrated in \eqref{LogCC1dstrong}. For that we need the following special class of binary encodings with adjacent elements (in the order induced by the SOS2 constraints) that only differ in one bit or coordinate. 
\begin{definition}\label{graycodecoro} We say $H=\bra{h^i}_{i=1}^n\in \mathcal{H}_{\lceil \log_2 n\rceil}(n)$ is a \emph{gray code} if and only if  for all $i\in \sidx{ n-1}$ we have $\sum\nolimits_{j=1}^{\lceil \log_2 n\rceil} \abs{h^{i}_j-h^{i+1}_j}=1$.
\end{definition}
The following corollary shows that formulation \eqref{SOS2char} for a gray code is precisely  the logarithmic formulation when $n$ is a power of two and otherwise eliminates some redundancy of the logarithmic formulation (see \cite{MuldoonPhd} for an alternate derivation and discussion about this redundancy). 
\begin{corollary}\label{sosub} Let $\mathcal{P}:=\bra{P^i}_{i=1}^n$ be the SOS2 constraint on $\Delta^{n+1}$ and  $H$ be a gray code, then  $\bra{\lambda,y}\in Q\bra{\mathcal{P},H}$ if and only if
\begin{subequations}\label{sos2formul}
\begin{alignat}{3}
\sum\nolimits_{j=1}^{n+1} \lambda_j &=1\label{sos2eqqq}\\
\sum\nolimits_{j=1}^{n+1} \min\set{h^j_l,h^{j-1}_l}\lambda_j\leq y_l &\leq \sum\nolimits_{j=1}^{n+1} \max\set{h^j_l,h^{j-1}_l}\lambda_j &\quad&\forall l \in \sidx{\left\lceil \log_2 n\right\rceil}\label{sos2formulgeneral}\\
\lambda_j&\geq 0 &\quad& \forall j\in \sidx{n+1}.\label{sos2formulbound}
\end{alignat}
\end{subequations}
where we let $h^0=h^1$ and $h^{n+1}=h^n$. Furthermore, all inequalities in \eqref{sos2formulgeneral} are facet defining and hence $\size_G\bra{Q\bra{\mathcal{P},H}}=2\left\lceil \log_2 n\right\rceil$. In addition, at least $n-1$ of the inequalities in \eqref{sos2formulbound} are facet defining and hence $\abs{\size_B\bra{Q\bra{\mathcal{P},H}} - n}\leq 1$ and  $\abs{\size\bra{Q\bra{\mathcal{P},H}}- \bra{2\left\lceil \log_2 n\right\rceil + n +2}}\leq 1$. Finally, these last two bounds are tight.  
\end{corollary} 
\proof{Proof.}
For any gray code  $c^i\in \bigcup_{j=1}^{\left\lceil \log_2 n\right\rceil}\set{-\e^j,\e^j}$ for each $i\in \sidx{n-1}$ and  $L(H)=\Real^{\left\lceil \log_2 n\right\rceil}$. Hence, the set of hyperplanes spanned by $\bra{c^i}_{i=1}^{n-1}$ is $\set{M\bra{b^l}}_{l=1}^{L}$ for $L=\left\lceil \log_2 n\right\rceil$ and $b^l= \e^i$ for each $l\in \sidx{L}$. This yields the results concerning non-bound inequalities.

For the bound on the number of facet defining inequalities in \eqref{sos2formulbound} note that  $L\bra{H}\neq\spann\bra{\set{c^i}_{i\in \sidx{n-1}\setminus\set{j-1} }}$ if and only if there exist $l\in \sidx{\left\lceil \log_2 n\right\rceil}$ such that $\abs{c^{j-1}}=\e^l$ and $\abs{c^{i}}\neq \e^l$ for all $i\in \sidx{n-1}\setminus\set{j-1}$. Then, the number of inequalities $\lambda_j\geq 0$ which fail to be facet defining is the number of bits that only change once in $H$, which is given by
\[T_n^1(H):=\abs{\set{l\in \sidx{\left\lceil \log_2 n\right\rceil}\,:\, \exists i_0\in \sidx{n-1} \text{ s.t. } h^{i_0}_l\neq h^{i_0+1}_l \text{ and } h^i_l=h^{i+1}_l \quad \forall l\in \sidx{n-1}\setminus \set{i_0}}}.\] 
To show that $T_n^1(H)\leq 2 $ assume for a contradiction that $n\geq 3$ and $T_n^1(H)\geq 3 $. Without loss of generality we may assume there exist $0=i_0< i_1<i_2<i_3<i_4=n$ such that for all $l\in \sidx{3}$ we have $h^{i}_l=1$ for all $i\in \sidx{i_l}$ and $h^i_{l}=0$ for all $i\in \sidx{i_l+1,n}$. Let $k=\lceil \log_2 n\rceil$ be the dimension or  number of bits of $H$ so that $n\geq 2^{k-1}+1$. Then, because the indices $i_1,i_2,i_3$ divide $H$ into $4=2^2$ sets, there exits $l\in\sidx{4}$ such that  $i_l-i_{l-1}\geq \left\lceil \bra{2^{k-1}+1}/4\right \rceil =2^{k-3}+1$. However, because $h^i\neq h^j$ and  $h^i_l=h^j_l$ for all $l\in \sidx{3}$ and $i,j\in \sidx{i_{l-1}+1,i_l}$, we must have $2^{k-3}\geq i_l-i_{l-1}\geq 2^{k-3}+1$ which yields a contradiction. For the tightness of the bounds note that for $H=\bra{(0,0)^T,(1,0)^T,(1,1)^T}$ we have $T_n^1(H)=2$ and for $H=\bra{(1,0,0)^T,(1,1,0)^T,(0,1,0)^T,(0,1,1)^T,(1,1,1)^T,(1,0,1)^T,(0,0,1)^T,(0,0,0)^T}$ we have $T_n^1(H)=0$. Finally, we recover the logarithmic formulation when $n$ is a power of two (e.g. formulation (9) in \cite{Modeling-Disjunctive-Constraints-FULL}) by subtracting \eqref{sos2eqqq} from the second inequality in \eqref{sos2formulgeneral}, noting that $ \min\set{h^j_l,h^{j-1}_l}=1$ if $h^j_l=h^{j-1}_l=1$ and $ \min\set{h^j_l,h^{j-1}_l}=0$ otherwise, and noting that $1-\max\set{h^j_l,h^{j-1}_l}=1$ if $h^j_l=h^{j-1}_l=0$ and $1-\max\set{h^j_l,h^{j-1}_l}=0$ otherwise. \Halmos\endproof

Finally, we can use Proposition~\ref{sos1prophyper} to give tight bounds on the embedding complexity of SOS2 constraints, which show that the logarithmic formulation from Corollary~\ref{sosub} is (nearly) optimal.
\begin{proposition}\label{sosboundcoro}Let $\mathcal{P}:=\bra{P^i}_{i=1}^n$ be the SOS2 constraint on $\Delta^{n+1}$, then
\begin{itemize}
\item $\mmc_G\bra{\mathcal{P}}=2\lceil \log_2 n\rceil$, and
\item  $n+3+\lceil \log_2 n\rceil\leq \mmc\bra{\mathcal{P}}\leq n+3 +2\lceil \log_2 n\rceil$.
\end{itemize}
\end{proposition}
\proof{Proof.}
Let $\set{b^l}_{l=1}^L\subseteq L(H)\setminus \set{\bf 0}$ be  such that $\set{M\bra{b^l}}_{l=1}^{L}$  is the set of linear hyperplanes spanned by $\set{c^i}_{i=1}^{n-1}$ in $L(H)$. By Proposition~\ref{sos1prophyper} an easy upper bound on the number of variable bounds $\lambda_j\geq 0$ that are not facet defining is $L$, which is achieved precisely if for all $l\in \sidx{L}$ there is only one $i\in \sidx{n-1}$ such that $b^l\cdot c^i\neq 0$. Hence we have  $\size\bra{Q\bra{\mathcal{P},H}}\geq 2+2L+n+1-L=n+3+L$. Then, to obtain a lower bound on both  $\mmc_G\bra{\mathcal{P}}$ and $\mmc\bra{\mathcal{P}}$ we need to minimize $L$. For that we note that $L$ is equal to the number of $1$-flats of the (central) hyperplane arrangement $\set{\set{b\in L(H)\,:\, c^i\cdot b =0}}_{i=1}^{n-1}$ in $L(H)$. Because $\spann\bra{\set{c^i}_{i=1}^{n-1}}=L(H)$ the number of such $1$-flats is at least $\dim(L(H))=\dim(H)$. Because all elements of $H$ are pairwise distinct  we have that $L\geq \dim(H)\geq \lceil \log_2 n\rceil$, which yields the lower bounds. The upper bounds follow from Corollary~\ref{sosub}.
\Halmos\endproof

The smallest size embedding formulation is achieved through the use of a very specific class of binary encoding. In the following section we explore how rare this class is with respect to yielding small embedding formulations.
  
\subsection{Size distribution for binary encodings}\label{graphsec}

If $n=2^k$ for some $k\in \Int$,  $\mathcal{P}:=\bra{P^i}_{i=1}^n$ is the SOS2 constraint and $H\in \mathcal{H}_{k}(n)$, then Proposition~\ref{sos1prophyper} shows that the  number of facets of $Q\bra{\mathcal{P},H}$ defined by general inequalities \eqref{sos2generalfacets} is upper bounded by $2 \binom{n-1}{k-1}$. The following proposition suggests that this upper bound may be nearly achieved. We include a proof of this result in Section~\ref{omproofapendix}.

\begin{restatable}{proposition}{antigraylemmaLem}\label{antigraylemma}Let $n=2^k$ for some $k\in \Int$ and $\mathcal{P}:=\bra{P^i}_{i=1}^n$ be the SOS2 constraint. There exist $H\in \mathcal{H}_{k}(n)$ such that $\size_G\bra{Q\bra{\mathcal{P},H}}$ is equal to twice the number of affine hyperplanes spanned by $\set{0,1}^{k-1}$.
\end{restatable}

\begin{figure}[htpb]
  \begin{center}
  \includegraphics[scale=0.37]{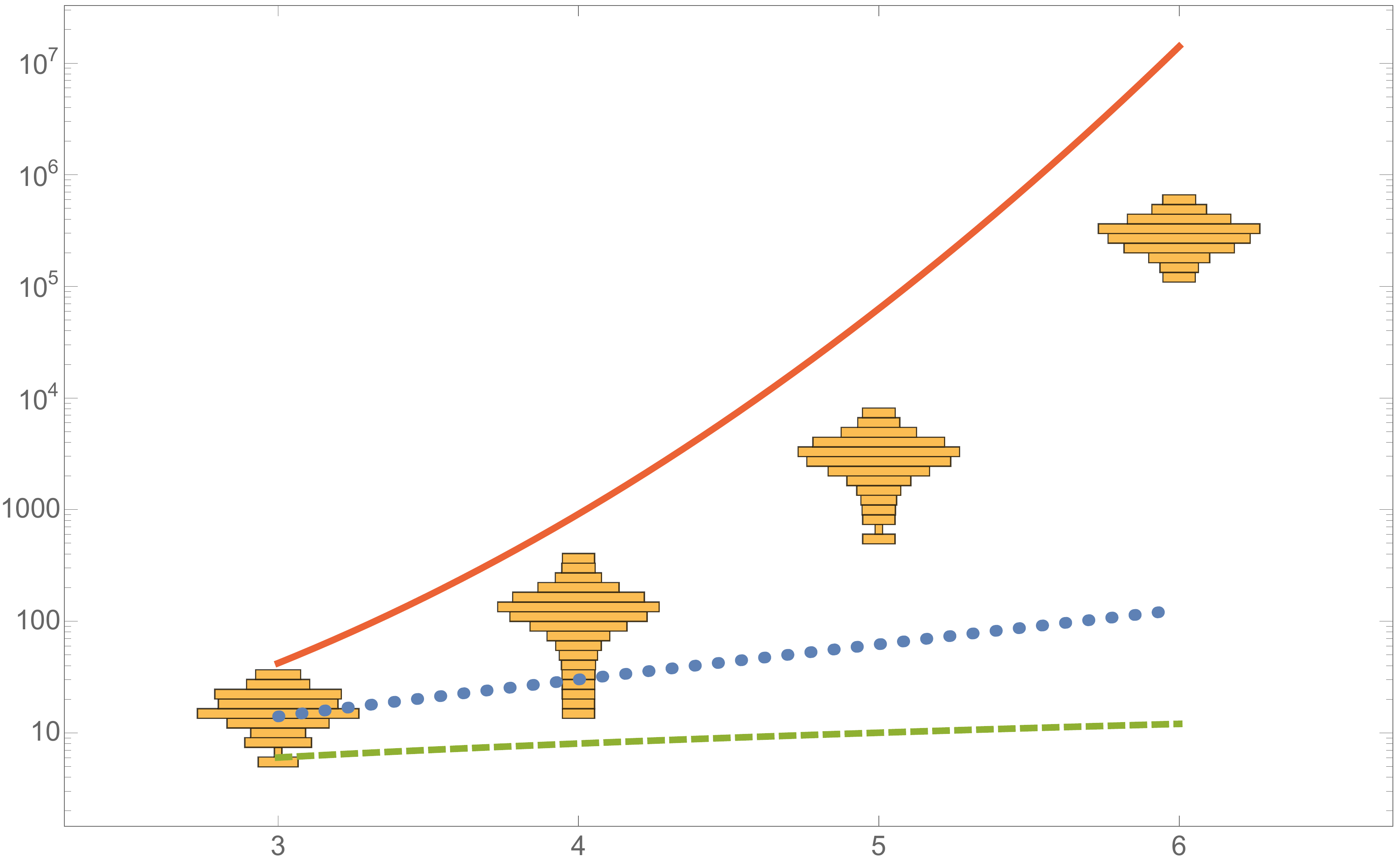}
  \end{center}
  \caption{$\size_G\bra{Q\bra{\mathcal{P},H}}$ for $k\in \sidx{3,6}$ and randomly selected binary encodings.}\label{randomfigure}
  \end{figure}

It is believed that the number of affine hyperplanes spanned by $\set{0,1}^{k-1}$ is close to its trivial upper bound of $\binom{n/2}{k-1}$ for $n=2^k$ (e.g. \cite{aichholzer1996classifying}). Both this upper bound and the general-inequality-defined facet bound of $2\binom{n-1}{k-1}$ grow roughly as $n^{\log_2n}$\footnote{We have $\Omega\bra{n^{(1-\varepsilon) \log_2 n}}=\binom{n/2}{k-1}\leq \binom{n-1}{k-1}\leq n^{\log_2 n}$ for all $\varepsilon>0$.}, which suggests that the worst case of $\size_G\bra{Q\bra{\mathcal{P},H}}$ for a binary encoded formulation is quasi-polynomial in $n$. 
Hence, it seems like an unfortunate selection of the specific binary encoding can lead to a formulation that is significantly larger than the lower bound from Proposition~\ref{sosboundcoro} or even the size of the unary encoded formulation from Corollary~\ref{unaryformulationsos2}. Because of its link with the number of hyperplanes spanned by subsets of $\set{0,1}^{k-1}$ (or $\set{-1,0,1}^{k-1}$), understanding the typical size of a binary encoded embedding formulation for SOS2 constraints may prove extremely challenging (e.g. \cite{voigt2006singular}). For this reason we only pursue a simple empirical study of the distribution of sizes for these formulations. For this study we selected $k\in \sidx{3,6}$ and calculated $\size_G\bra{Q\bra{\mathcal{P},H}}$ for randomly selected binary encodings (the ones associated to a random permutation of $\set{0,1}^{k}$). For $k=3$ we considered all $40,320$ possible encodings, while for $k\in \set{4,5}$ we only used a random sample of $10,000$ encodings and for $k=6$ we only used a random sample of $1,000$ encodings (calculating the formulation sizes for $k=6$ was already computational intensive). The results of this study are presented in Figure~\ref{randomfigure}.
 The figure presents histograms for $\size_G\bra{Q\bra{\mathcal{P},H}}$ with random binary encoding $H$ for each $k$, together with the trivial upper bound of $2 \binom{n-1}{k-1}$ (depicted by the solid red line),  $\size_G\bra{Q\bra{\mathcal{P},H}}$ for the unary encoding (depicted by the dotted blue line) and $\size_G\bra{Q\bra{\mathcal{P},H}}$ for the optimal binary encoded formulation (depicted by the dashed green line). The figure shows that the typical value of $\size_G\bra{Q\bra{\mathcal{P},H}}$ for a binary encoding seems to be much closer to the upper bound and suggests that a randomly selected encoding may often lead to a formulation that is significantly larger than even the unary encoded formulation. Hence a careful encoding selection appears crucial to obtain a small formulation.

\section{Formulations for Piecewise Linear Functions of Two Variables}\label{pwlsecc}

\begin{figure}[htpb]
  \begin{center}
  \subfigure[Piecewise Linear Function on the Union-Jack Triangulation for $m=8$]{\label{newpwl2fig1}\includegraphics[scale=0.3]{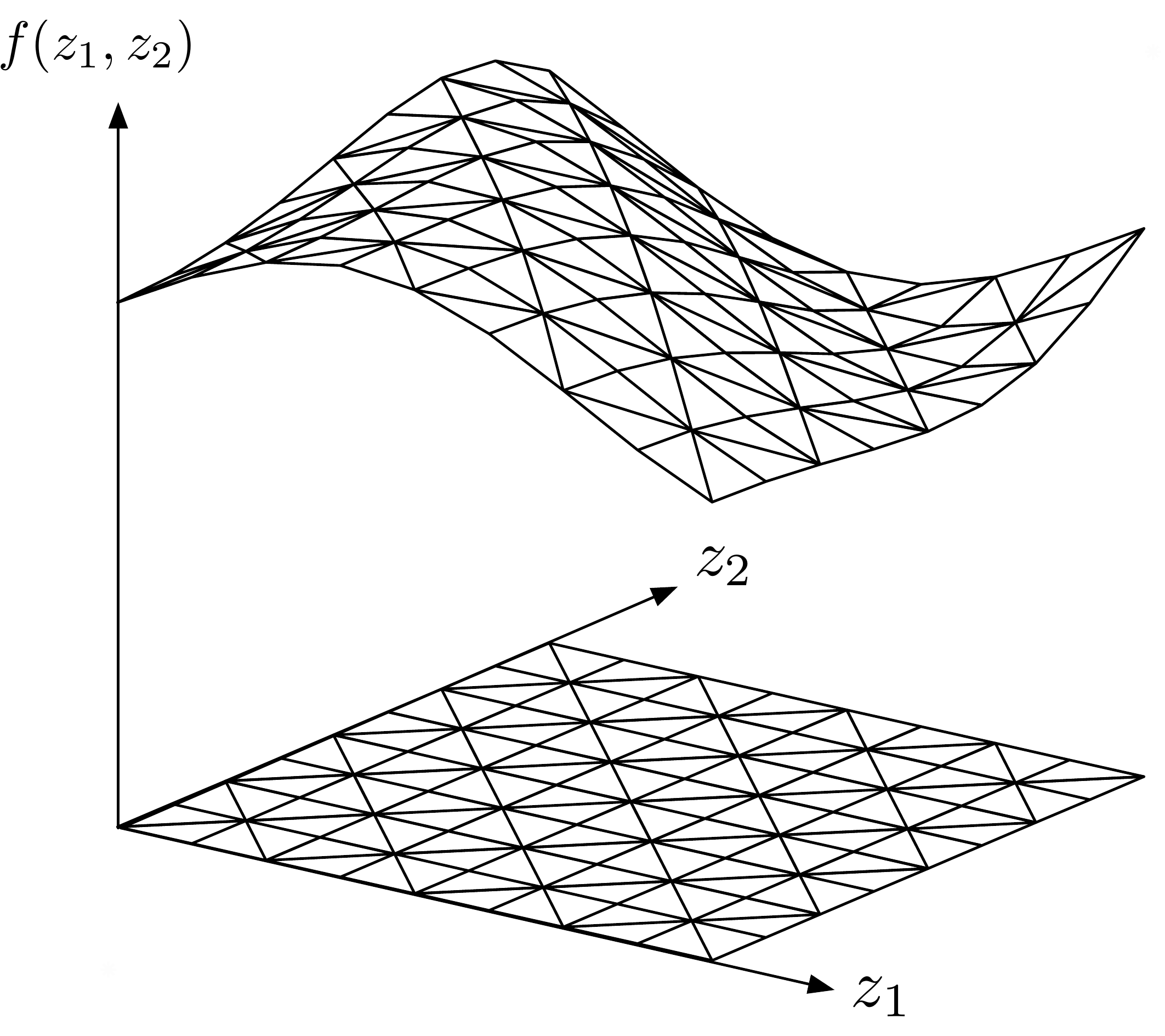}}\quad\quad\quad 
    \subfigure[The Union-Jack Triangulation for $m=2$.]{\label{newpwl2fig2}\includegraphics[scale=0.2]{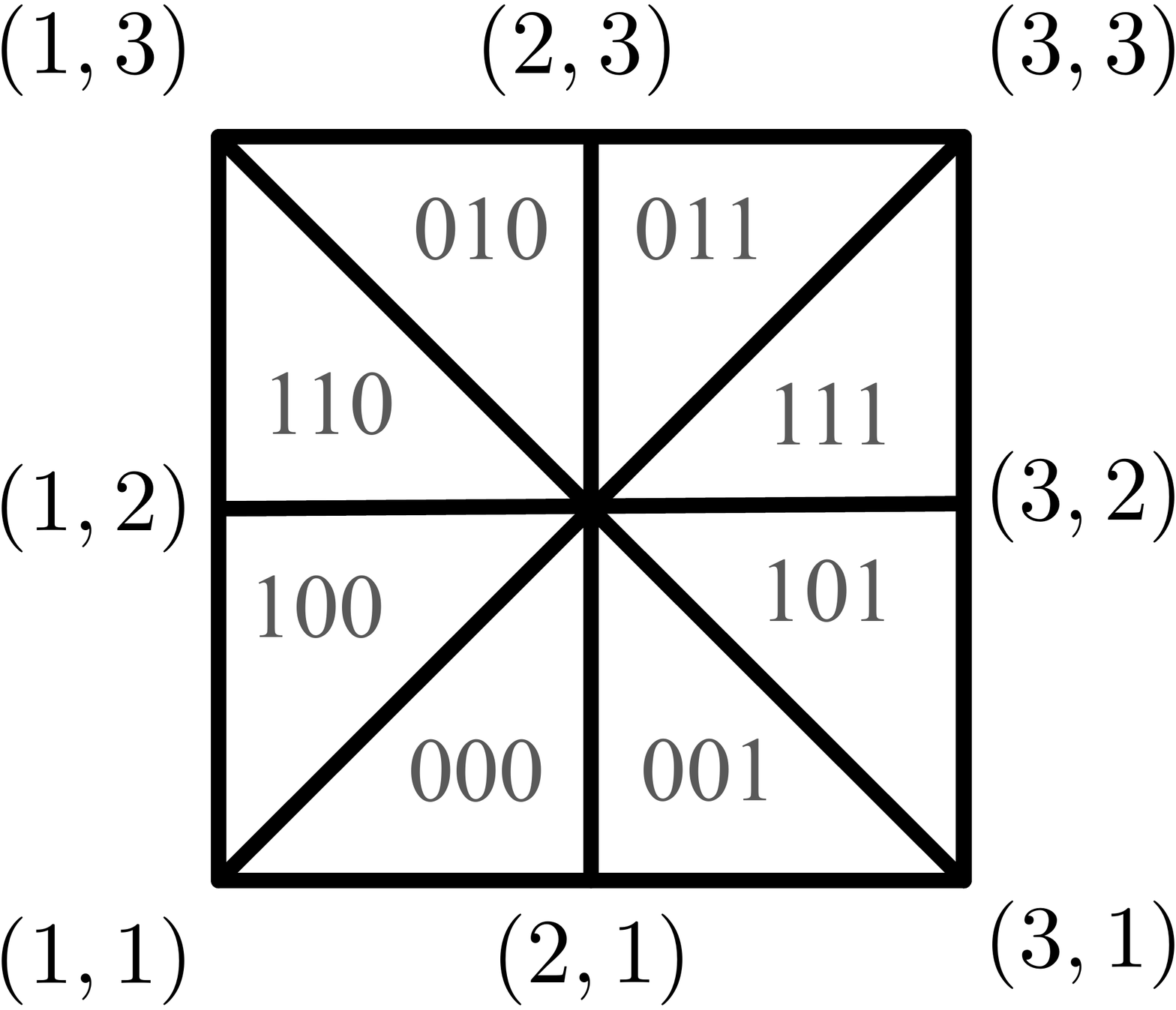}}
  \end{center}
  \caption{Piecewise Linear Functions of Two Variables and Grid Triangulations.}\label{newpwl2fig}
  \end{figure}

The results in Section~\ref{graphsec} show that it may be hard to construct small embedding formulations. However, we now show how small embedding formulations can be constructed for multivariate piecewise linear functions. MIP formulations for multivariate piecewise linear functions can be constructed using standard generalizations of SOS2 constraints (e.g.  \cite{lee01,Modeling-Disjunctive-Constraints-FULL,Mixed-Integer-Models-for-Nonseparable}). For simplicity,  we only consider formulations for piecewise linear functions of two variables defined on grid triangulations on $\sidx{m+1}^2$ such as those depicted in Figure~\ref{newpwl2fig}. More precisely, we consider functions $f:[1,m+1]^2\to \Real$ that are continuous in $[1,m+1]^2$ and affine in each triangle of the triangulation (e.g. for the triangulation in Figure~\ref{newpwl2fig2} it is affine in $\conv\bra{\set{\bra{1,1},\bra{1,2},\bra{2,2}}}$). 

The following proposition summarizes the standard generalization of SOS2 constraints used to give a disjunctive representation of the graph of a piecewise linear function. It also describes how an embedding formulation can be used to model the corresponding disjunctive constraint and yield a formulation for the graph of this function.

\begin{proposition}\label{PWLFORMPROP}Let $\mathcal{T}=\set{T^1_{u,v},T^2_{u,v}}_{u,v=1}^{m}$ be a grid triangulation on $\sidx{m+1}^2$ so that for all $u,v\in \sidx{m}$ we have $T^1_{u,v}\cup T^2_{u,v}= \set{u,u+1}\times\set{v,v+1}$, $\abs{T^1_{u,v}}=\abs{T^2_{u,v}}=3$ and either $T^1_{u,v}\cap T^2_{u,v}=\set{(u,v),(u+1,v+1)}$ or  $T^1_{u,v}\cap T^2_{u,v}=\set{(u+1,v),(u,v+1)}$ (e.g. for the triangulation in Figure~\ref{newpwl2fig2} we can take $T_{2,1}^1=\set{\bra{2,1},\bra{3,1},\bra{2,2}}$ and $T_{2,1}^2=\set{\bra{3,1},\bra{3,2},\bra{2,2}}$ or exchange $T_{2,1}^1$ and $T_{2,1}^2$). Finally, let $\Delta_2^{m+1}:=\set{\lambda\in \Real^{\sidx{m+1}^2}_+\,:\, \sum_{u,v=1}^{m+1} \lambda_{\bra{u,v}}=1}$, $n=2m^2$ and  $\mathcal{P}\bra{\mathcal{T}}:=\bra{P^i}_{i=1}^{n}$ so that $P^{2\bra{m (u-1)+(v-1)}+t}=P\bra{T_{u,v}^t}:=\set{\lambda\in \Delta_2^{m+1}\,:\,  \lambda_{\bra{\bar{u},\bar{v}}}\leq 0\quad \forall \bra{\bar{u},\bar{v}}\notin T_{u,v}^t}$ for all $u,v\in\sidx{m}$ and $t\in \set{1,2}$. 

Then for any continuous function $f:[1,m+1]\to \Real$ that is affine in $\conv\bra{T}$ for each $T\in \mathcal{T}$, we have that a disjunctive representation of its graph $\gr(f):=\set{\bra{x,z}\in \Real^3\,:\, f(x)=z}$ is given by  
  \begin{subequations}\label{vformexunionpwl}
      \begin{alignat}{3}
    \sum\nolimits_{u,v=1}^{m+1}  u \lambda_{\bra{u,v}}=x_1,\quad  \sum\nolimits_{u,v=1}^{m+1}  v \lambda_{\bra{u,v}}=x_2,\quad \sum\nolimits_{u,v=1}^{m+1}  f(u,v) \lambda_{\bra{u,v}}&=z&\label{vformnumpwl}\\
            \lambda&\in \bigcup\nolimits_{i=1}^{n} P^i.\label{vformnumpwldisj}
    \end{alignat}
  \end{subequations}

If $k\geq \left\lceil \log_2 n\right\rceil$ and $H\in \mathcal{H}_k(n)$ then an ideal formulation of $\gr(f)$ is given by \eqref{vformnumpwl}, $\bra{\lambda,y}\in Q\bra{\mathcal{P}\bra{\mathcal{T}},H} $ and $y\in \mathbb{Z}^k$. 
\end{proposition}

Constructing and analyzing  embedding formulations for triangulations can be significantly more complicated than for SOS2 constraints. However, the analysis can be partially achieved through the following simple lemma that shows how a description of  $Q\bra{\mathcal{P}\bra{\mathcal{T}},H}$ can be obtained from an ideal formulation of \eqref{vformnumpwldisj} whose binary variables are \emph{compatible} with the encoding $H$. That is if the formulation can be re-interpreted as a formulation of $(\lambda,y)\in \bigcup\nolimits_{i=1}^{n} P_i \times \set{h^i}$ (the LP relaxation of such formulation corresponds to $Q$ in the lemma).
\begin{lemma}\label{formulationLemma}Let $\mathcal{P}:=\bra{P^i}_{i=1}^n$ be a finite family of polyhedra in $\Real^d$ satisfying Assumption~\ref{ass1}, $k\geq \left\lceil \log_2 n\right\rceil$,  
$H:=\bra{h^i}_{i=1}^n\in \mathcal{H}_k(n)$ and $Q\subseteq \Real^{d+k}$ be a rational polyhedron. If $\ext\bra{Q}\subseteq \Real^d\times \Int^k$ and
\begin{equation}\label{lemmacondition}\bra{x,y}\in Q\cap \bra{\Real^d\times\Int^k} \quad \Leftrightarrow \quad\exists i\in \sidx{n} \text{ s.t. } y=h^i\;\wedge\; x\in P^i,\end{equation}
then $Q=Q\bra{\mathcal{P},H}$.
\end{lemma}
\proof{Proof.} By \eqref{lemmacondition}, the definition of $Q\bra{\mathcal{P},H}$ and convexity of $Q$ we have $Q\bra{\mathcal{P},H}\subseteq Q$. Furthermore, \eqref{lemmacondition} and the assumption on $\ext\bra{Q}$ yield $\ext\bra{Q}\subseteq \bigcup_{i=1}^n P^i\times \set{h^i}$. By Proposition~\ref{firstprop} and the fundamental theorem of integer programming applied to $Q\cap \bra{\Real^d\times\Int^k}$ and $Q\bra{\mathcal{P},H}\cap \bra{\Real^d\times\Int^k}$  we further have $Q_\infty=P^1_\infty\times\set{\mathbf{0}}=Q\bra{\mathcal{P},H}_\infty$. Finally, combining these properties with Minkowski-Weyl we get $Q=\conv\bra{\ext\bra{Q}}+Q_\infty\subseteq Q\bra{\mathcal{P},H}$, which shows the result.
\Halmos\endproof
Combining Lemma~\ref{formulationLemma} and formulations from  \cite{lee01} and \cite{Modeling-Disjunctive-Constraints-FULL} we can give a rather precise analysis for a triangulation known as the \emph{union-jack} \citep{todd77}, which is depicted in Figure~\ref{newpwl2fig2} (See \cite{Modeling-Disjunctive-Constraints-FULL} for a precise description). In particular, we have that using the unary encoding leads to an extremely large embedding formulation, but using a carefully selected binary encoding leads to a simple formulation with near-optimal size. In the following proposition we let $\size_B\bra{Q\bra{\mathcal{P},H}}$ denote the number of facets of $Q\bra{\mathcal{P},H}	$ defined by variable bounds as we did for SOS2 constraints.

\begin{restatable}{proposition}{unionjackboundprop}\label{UnionJackCoro}Let $\mathcal{P}=\mathcal{P}(\mathcal{T})$ for a grid triangulation $\mathcal{T}$ on $\sidx{m+1}^2$ and $n=2m^2$ so that $\abs{\mathcal{P}\bra{\mathcal{T}}}=n$. Then $\bra{\sqrt{n/2}+1}^2\leq \mmc\bra{\mathcal{P}}$ and if $H$ is the unary encoding then
\[\size_B\bra{Q\bra{\mathcal{P},H}}=\bra{\sqrt{n/2}+1}^2  \quad\text{and}\quad \size\bra{Q\bra{\mathcal{P},H}}=4+\bra{\sqrt{n/2}+1}^2+ \binom{2\sqrt{n/2}}{\sqrt{n/2}}.\]
In contrast, if $m$ is a power of two and $\mathcal{T}$ is the union-jack triangulation, then there exist a binary encoding $H\in \mathcal{H}_{\log_2 n}(n)$ such that 
\begin{equation}\label{unionjackopt}
\size_B\bra{Q\bra{\mathcal{P},H}} = \bra{\sqrt{n/2}+1}^2\quad\text{and}\quad \size\bra{Q\bra{\mathcal{P},H}} = 4+\bra{\sqrt{n/2}+1}^2+ 2  \log_2 (n/2).
\end{equation}
\end{restatable}

We postpone a formal proof of Proposition~\ref{UnionJackCoro} to Section~\ref{proofofUnionJackCoro} and instead illustrate it with the following example that shows how Lemma~\ref{formulationLemma} can be used to recover an embedding formulation from the ideal formulation from \cite{Modeling-Disjunctive-Constraints-FULL}. In particular, it shows how studying this formulation reveals the encoding needed to recover the embedding formulation.
\begin{example}\label{unionjackex} The union-jack triangulation $\mathcal{T}=\set{T^1_{u,v},T^2_{u,v}}_{u,v=1}^{m}$ for $m=2$ depicted in Figure~\ref{newpwl2fig2} can be described in a standard format used to construct the formulation in \cite{Modeling-Disjunctive-Constraints-FULL} by letting $T^1_{u,v}=\set{(2,2),(2,2v-3),(2u-3,2v-3)}$ and  $T^2_{u,v}=\set{(2,2),(2u-3,2),(2u-3,2v-3)}$ for all $u,v\in \sidx{2}$. In addition, the formulation from \cite{Modeling-Disjunctive-Constraints-FULL} for this triangulation is given by 
 \begin{subequations}\label{unionjackexform}
 \begin{alignat}{6}
   \lambda_{\bra{2,1}}+\lambda_{\bra{2,3}}&\leq 1-y_1,&\quad 
  \lambda_{\bra{1,2}}+\lambda_{\bra{3,2}}&\leq y_1&\quad& \\
    \lambda_{\bra{1,1}}+\lambda_{\bra{2,1}}+\lambda_{\bra{3,1}}& \leq 1-y_2, &\quad 
\lambda_{\bra{1,3}}+\lambda_{\bra{2,3}}+\lambda_{\bra{3,3}}& \leq y_2 &\quad&\\
  \lambda_{\bra{1,1}}+\lambda_{\bra{1,2}}+\lambda_{\bra{1,3}}& \leq 1-y_3, &\quad  
\lambda_{\bra{3,1}}+\lambda_{\bra{3,2}}+\lambda_{\bra{3,3}}& \leq y_3 &\quad&\\
 y &\in \set{0,1}^3,&\quad \sum\nolimits_{u,v=1}^3 \lambda_{(u,v)} &=1, &\quad& \lambda_{(u,v)} &\geq 0 \quad \forall u,v\in \sidx{3}.
 \end{alignat}
 \end{subequations}
If we let $\mathcal{P}\bra{\mathcal{T}}:=\bra{P^i}_{i=1}^{8}$ as defined in Proposition~\ref{PWLFORMPROP}, we can check that if $\bra{y,\lambda}$ is feasible for \eqref{unionjackexform} and $y=(0,0,0)$ then $\lambda \in P^1=P\bra{T^1_{1,1}}=\set{\lambda\in \Delta_2^3\,:\, \lambda_{\bra{u,v}}\leq 0\quad \forall \bra{u,v}\notin \set{\bra{2,2},\bra{2,1} ,\bra{1,1}}}$. Then letting $h^1=(0,0,0)$ we obtain condition \eqref{lemmacondition} of Lemma~\ref{formulationLemma} for $i=1$. Similarly, we may iterate over all  values of $y \in \set{0,1}^3$ to obtain the complete  $P\bra{T^t_{u,v}}$-$h^{2\bra{m (u-1)+(v-1)}+t}$ or triangle-vector assignment for Lemma~\ref{formulationLemma}  depicted in Figure~\ref{newpwl2fig2} and given by  
 \begin{alignat*}{3}
 h^1=(0,0,0),\; T_{1,1}^1=\set{\bra{2,2},\bra{2,1} ,\bra{1,1}};&\quad\quad\quad&h^2=(1,0,0),\;  T_{1,1}^2=\set{\bra{2,2},\bra{1,2},\bra{1,1}};\\
 h^3=(0,0,1),\; T_{2,1}^1=\set{\bra{2,2},\bra{2,1},\bra{3,1}};&\quad\quad\quad&h^4=(1,0,1),\;  T_{2,1}^2=\set{\bra{2,2},\bra{3,2},\bra{3,1}};\\ 
 h^5=(0,1,0),\; T_{1,2}^1=\set{\bra{2,2},\bra{2,3},\bra{1,3}};&\quad\quad\quad&h^6=(1,1,0),\; T_{1,2}^1=\set{\bra{2,2},\bra{1,2},\bra{1,3}};\\ 
 h^7=(0,1,1),\; T_{2,2}^1=\set{\bra{2,2},\bra{2,3},\bra{3,3}};&\quad\quad\quad&h^8=(1,1,1),\; T_{2,2}^1=\set{\bra{2,2},\bra{3,2},\bra{3,3}}.
  \end{alignat*}
  Because formulation \eqref{unionjackexform} is ideal, by Lemma~\ref{formulationLemma} we have that its LP relaxation is equal to $Q\bra{\mathcal{P}\bra{\mathcal{T}},H}$ for this $H$.
\end{example}

\begin{figure}[htpb]
  \begin{center}
  \subfigure[Modified Union-Jack Triangulation for $m=4$.]{\label{triag2a}\quad\quad\quad\quad\quad\quad\includegraphics[scale=.18]{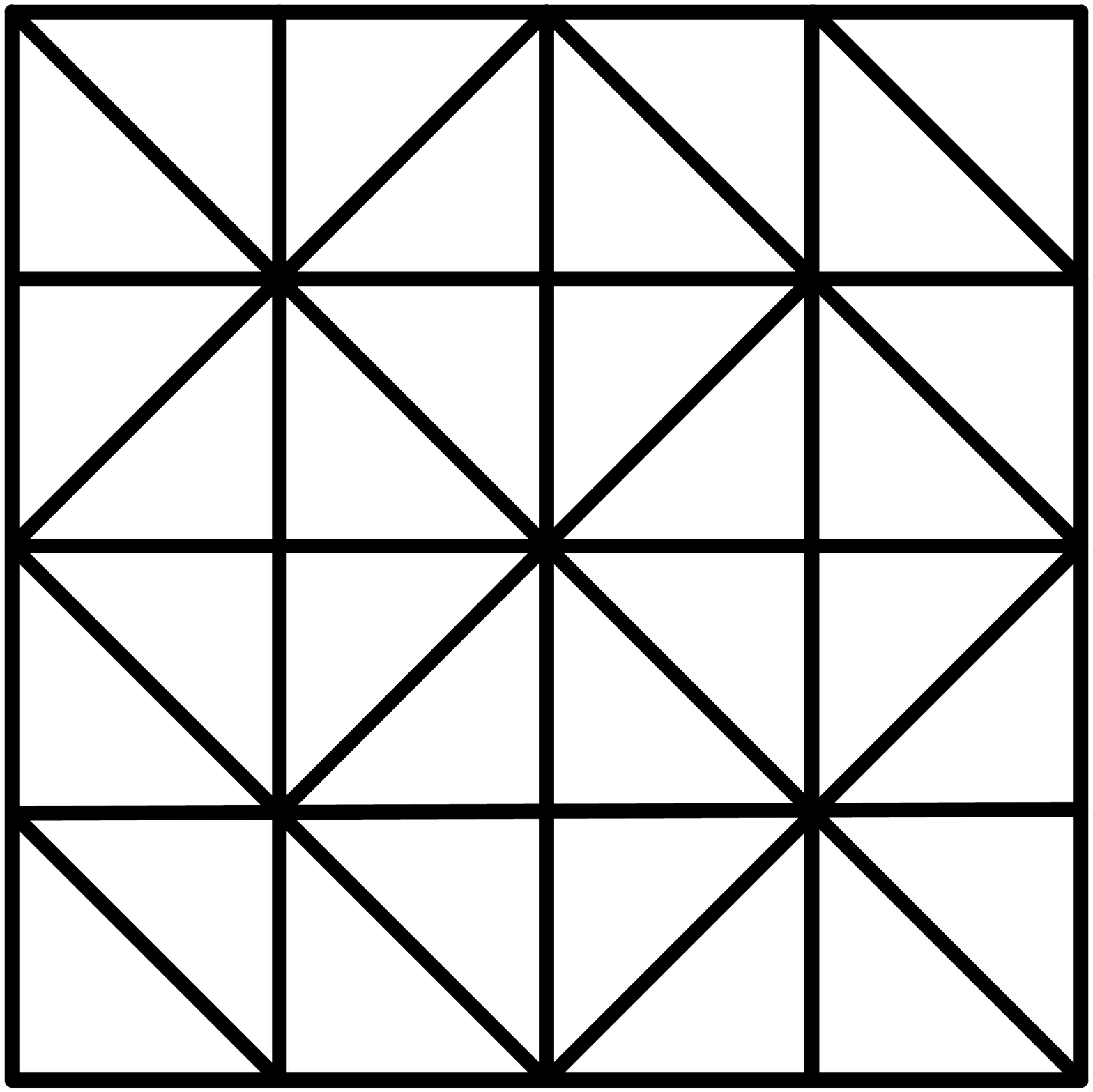}\quad\quad\quad}
    \subfigure[K1 Triangulation for for $m=4$.]{\label{triag2b}\quad\quad\quad\includegraphics[scale=.18]{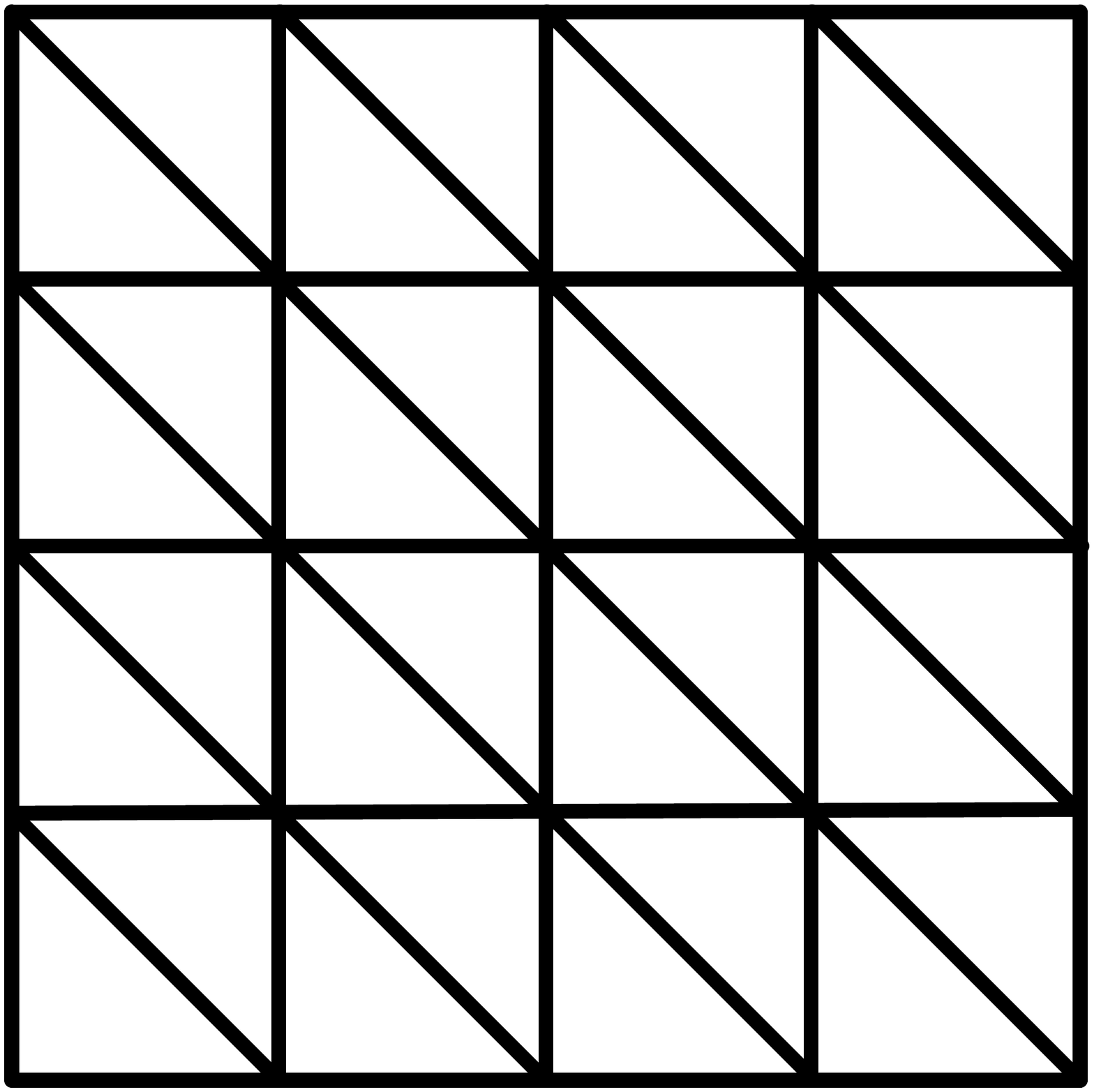}\quad\quad\quad\quad\quad\quad}
  \end{center}
  \caption{Different Triangulations.}\label{triag2}
  \end{figure}

 Proposition~\ref{UnionJackCoro} shows that the specific encoding used can have a significant impact on the size of an embedding formulation for triangulations. As illustrated in Example~\ref{unionjackex}, we can use an existing ideal formulation and Lemma~\ref{formulationLemma} to recover a favorable encoding.   Unfortunately, the formulation from \cite{Modeling-Disjunctive-Constraints-FULL} used to obtain the favorable encoding only works for the union-jack triangulation and it is sometimes preferable to use different triangulations such as the ones depicted in Figure~\ref{triag2} (e.g. \cite{Fitting-Piecewise-Linear}). In the following subsection we explore how adapting the favorable encoding for the union-jack triangulation to similar triangulations can sometimes help computationally construct a small embedding formulation.

\subsection{Constructing Embedding Formulations Computationally}\label{cddsec}

One way to construct embedding formulations is to computationally construct the convex hull in \eqref{embconv} for a specific encoding. Picking a random encoding will likely result in an extremely large formulation (cf. Section~\ref{graphsec}). For this reason we now investigate the effectiveness of using a known favorable encoding for a similar constraint. For this we consider the modification of the union-jack triangulation illustrated in Figure~\ref{triag2a}. This triangulation is obtained by changing the way the bottom-left and top-right squares of the triangulation are divided into two triangles (in the original triangulation they are divided into top-left and bottom-right triangles, and in the modified triangulation they are divided into bottom-left and top-right triangles). 

\begin{figure}[htpb]
  \begin{center}
    \subfigure[Encoding Adaptation for $m=2$.]{\label{aaa}\includegraphics[scale=.38,trim=0 -1.2cm 0 0]{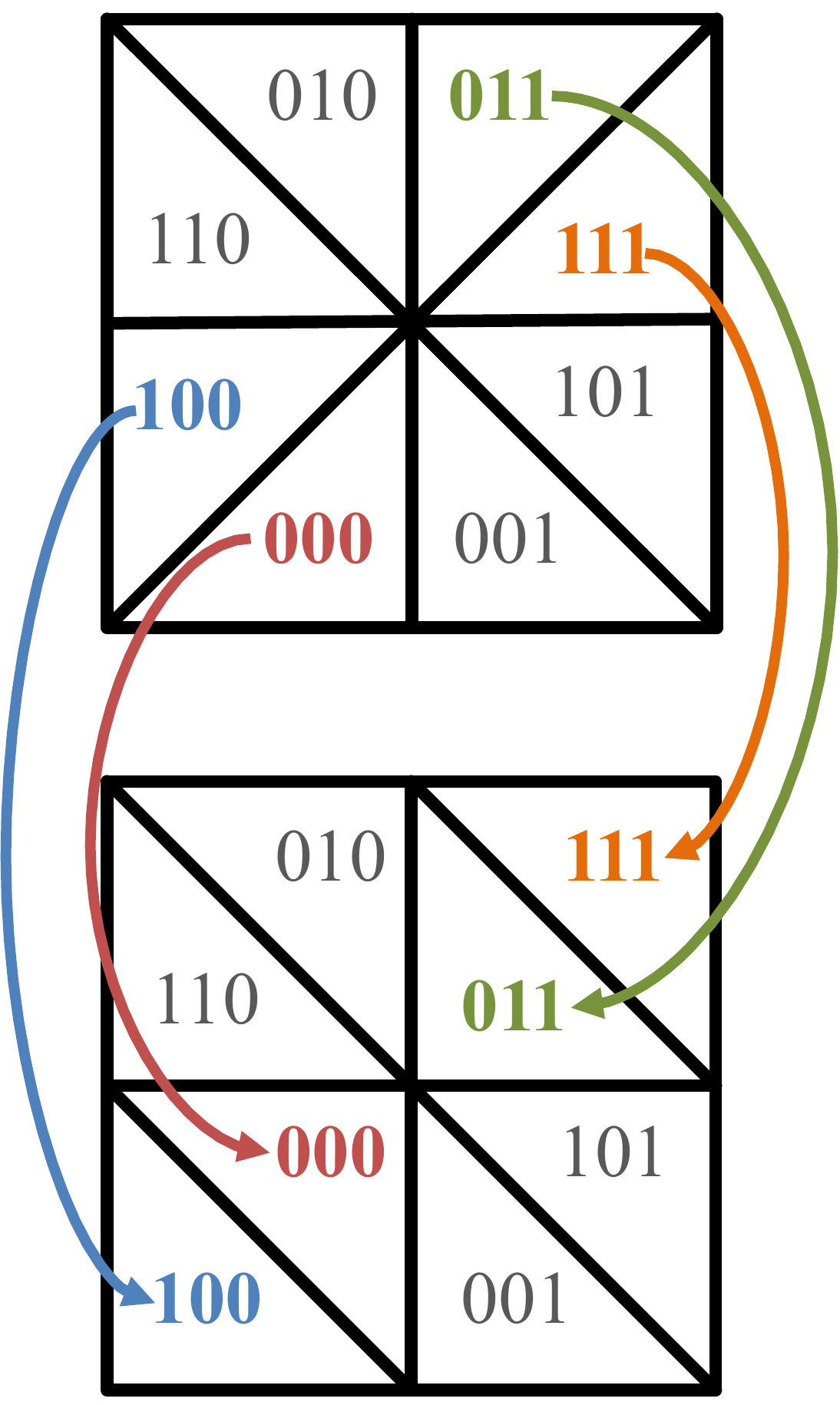}\quad\quad\quad}
    \subfigure[Solve Times for $m\in \set{4,8,16,32}$ {[s]}]{\label{resfig1}\quad\quad\includegraphics[scale=1.1]{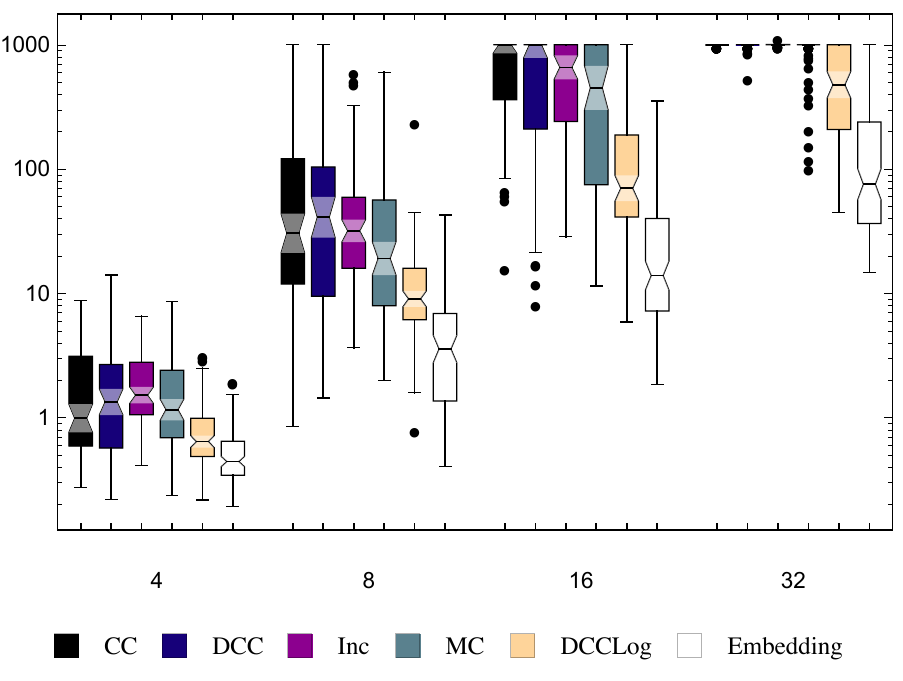}}
  \end{center}
  \caption{Encoding and Solve Times for Modified Union-Jack Triangulation.}
  \end{figure}

To construct an embedding formulation for this modified triangulation we adapt the encoding $H$ associated to the logarithmic formulation for the union-jack triangulation (i.e. the encoding illustrated in Example~\ref{unionjackex}). The adaptation uses the same triangle-vector assignment for all triangles except for the ones belonging to the bottom-left and top-right squares. This adaptation is illustrated in Figure~\ref{aaa} for $m=2$ where the top triangulation is the union-jack with the encoding described in Example~\ref{unionjackex} and the bottom triangulation is the modified  union-jack with the adapted encoding. For both squares where the encoding is changed, the adaptation assigns to the bottom-left triangle the same $h^i$ assigned to the top-left triangle in the original encoding for the union-jack triangulation. Similarly, the adaptation assigns to top-right triangle the same $h^i$ assigned to the bottom-right triangle in the original encoding. To construct the formulation we simply compute the convex hull of $\bigcup_{=1}^n P^i \times \set{h^i}$ for this modified $H=\bra{h^i}_{i=1}^n$  and $\mathcal{P}\bra{\mathcal{T}}:=\bra{P^i}_{i=1}^{n}$ from Proposition~\ref{PWLFORMPROP} using the software cddlib \citep{cddlib}.

We tried this for $m\in \set{4,8,16,32}$ and for all four cases the resulting embedding formulation only had four more inequalities than the formulation for the original union-jack triangulation. In addition, computing the convex hulls with cddlib for each $m\in \set{4,8,16,32}$ took respectively less than a second, 10 seconds, 24 minutes and 3.5 days on an Intel i7-3770 3.40GHz workstation with 32GB of RAM. The computational time can grow quickly with $m$, but fortunately this computation only has to be done once and the formulation can then be stored. Similar to traditional MIP formulations, the resulting stored embedding formulation can be used for free in any problem that requires a piecewise linear functions of two variables based on the modified union-jack triangulation. More specifically, the same formulation can be used independently of the specific data (e.g. actual function values) associated to the piecewise linear functions (i.e. in the formulation from Proposition~\ref{PWLFORMPROP} we need to update \eqref{vformnumpwl}, but we do not need to change the embedding formulation for disjunction \eqref{vformnumpwldisj}).  Hence, 3.5 days does not seem that large when compared with the research time required to develop a small and ideal  ad-hoc MIP formulation. Of course, this statement is conditional on the resulting formulation being small (so that it can be effectively stored and reused) and the formulation yielding a computational advantage (which is correlated with, but not guaranteed by a small size).  

To check if the formulation obtained with cddlib preserves the computational advantage of the logarithmic formulation for the original  union-jack triangulation, we replicate the computational experiments in   \cite{Modeling-Disjunctive-Constraints-FULL} and \cite{Mixed-Integer-Models-for-Nonseparable} for the modified union-jack triangulation. These experiments consider a series of transportation problems whose objective functions are the sum of $25$ piecewise linear functions of two variables on the $\sidx{m+1}^2$ grid for $m\in \set{4,8,16,32}$. For each $m$ the experiment considers $100$ randomly generated instances. With the exception of the logarithmic formulation, all formulations considered  in the original experiment are applicable for the modified union-jack triangulation. So we test all these formulations with the logarithmic formulation replaced by the embedding formulation constructed computationally using cddlib. All formulations were implemented using the JuMP modeling language \citep{jump,LubinDunningIJOC,DunningHuchetteLubin2015} and solved with Gurobi v6.5 \citep{gurobi} on an Intel i7-3770 3.40GHz workstation with 32GB of RAM. Solve times for all combinations of formulations and solver are presented in Figure~\ref{resfig1} for $m\in \set{4,8,16,32}$. We refer the reader to \cite{Mixed-Integer-Models-for-Nonseparable} for details on the benchmark formulations, but we note that DCC and MC is obtained by variants of Theorem~\ref{BTheo}, CC is the generalization of formulation \eqref{CC1d} and DCCLog is obtained by combining a variant of Theorem~\ref{BTheo} with the same encoding used for the embedding formulation. We can see that the embedding formulation can provide a significant computational advantage.
\section{Omitted Proofs and Additional Examples}

\subsection{Proof of Proposition~\ref{sos1prophyper}}\label{proofofsos1prophyper}
To prove Proposition~\ref{sos1prophyper} we use the following lemma, which gives a precise characterization of the facial structure of $Q\bra{\mathcal{P},H}$ for SOS2 constraints.
\begin{lemma}\label{facelemma} Let $\mathcal{P}:=\bra{P^i}_{i=1}^n$ be the SOS2 constraint on $\Delta^{n+1}$, $H:=\bra{h^i}_{i=1}^n\in \mathcal{H}_k(n)$, $h^0=h^1$, $h^{n+1}=h^n$, $c^i=h^{i+1}-h^i$ for $i\in \sidx{0,n}$. In addition, for any $J^-,J^+\subseteq \sidx{n+1}$ let 
\begin{alignat*}{3}
E\bra{J^-,J^+}:&=\bigcup\nolimits_{j\in J^-}\set{\bra{\e^j,h^{j-1}}} &\cup \bigcup\nolimits_{j\in J^+}\set{\bra{\e^j,h^{j}}}\\ &=\bigcup\nolimits_{j\in J^-}\set{\bra{\e^j,h^{j}-c^{j-1}}} &\cup \bigcup\nolimits_{j\in J^+}\set{\bra{\e^j,h^{j}}}
\end{alignat*}
so that $Q\bra{\mathcal{P},H}=\conv\bra{E\bra{\sidx{n+1},\sidx{n+1}}}$ and $E\bra{\sidx{n+1},\sidx{n+1}}=\ext\bra{Q\bra{\mathcal{P},H}}$.

Then for any for any $J^-,J^+\subseteq \sidx{n+1}$ we have  \[\dim\bra{E\bra{J^-,J^+}}=\abs{J^+\cup J^-}-1+\dim\bra{\set{c^{j-1}}_{j\in J^+\cap J^-}}.\]  Furthermore, $F\subseteq \Real^{n+1+k}$ is a face of $Q\bra{\mathcal{P},H}$ if and only if there exist  $J^-,J^+\subseteq \sidx{n+1}$ such that $F=F\bra{J^-,J^+}:=\conv\bra{E\bra{J^-,J^+}}$ and there exist $b\in \Real^k$ such that 
\begin{equation}\label{facelemmacond}
b\cdot c^{j-1}=0 \quad \forall j\in J^+\cap J^-,\quad b\cdot c^{j-1}< 0 \quad j\in J^-\setminus J^+ \quad\text{and}\quad b\cdot c^{j-1}> 0\quad \forall j\in J^+\setminus J^-.
\end{equation}
\end{lemma}
\proof{Proof.}
For the dimension  of $E\bra{J^-,J^+}$, let  $J\subseteq J^-\cap J^+$ be such that $\spann\bra{\set{c^{j-1}}_{j\in J}}=\spann\bra{\set{c^{j-1}}_{j\in J^-\cap J^+}}$ and $\abs{J}=\dim\bra{\set{c^{j-1}}_{j\in J^-\cap J^+}}$, and 
\[E=\bigcup\nolimits_{j\in J^-\setminus J^+}\set{\bra{\e^j,h^{j-1}}} \cup \bigcup\nolimits_{j\in J^+\setminus J^-}\set{\bra{\e^j,h^{j}}}\cup \bigcup\nolimits_{j\in J^-\cap J^+}\set{\bra{\e^j,h^{j}}}\cup \bigcup\nolimits_{j\in J}\set{\bra{\e^j,h^{j-1}}}. \]
We can check that $E$ is a set of $\abs{J^+\cup J^-}+\dim\bra{\set{c^{j-1}}_{j\in J^+\cap J^-}}$ affinely independent vectors. Finally, for any $\bra{\lambda,y}\in E\bra{J^-,J^+}\setminus E$ there exist $j\in J^+\cap J^-\setminus J$ such that $\bra{\lambda,y}=\bra{\e^j,h^{j-1}}=\bra{\e^j,h^{j}-c^{j-1}}$. Let $\mu\in \Real^{J}$ be such that $c^{j-1}=\sum_{i\in J} \mu_i c^{i-1}$. Then $\bra{\e^j,h^{j}-c^{j-1}}=
\mu_0 \bra{\e^j,h^{j}}+\sum_{i\in J} \bra{ \mu_i \bra{\e^i,h^i- c^{i-1}} -\mu_i \bra{\e^i,h^i} }$ for $\mu_0=1$ and result follows. 

For the facial characterization we have that $F\subseteq \Real^{n+1+k}$ is a face of $Q\bra{\mathcal{P},H}$ if and only if there exist $\bra{a,b}\in \Real^{n+1+k}$ such that 
\begin{equation}\label{facemax}
\begin{aligned}F&=\argmax\set{a\cdot \lambda+b\cdot y\,:\, \bra{\lambda,y}\in Q\bra{\mathcal{P},H}}\\&=\conv\bra{\argmax\set{a\cdot \lambda+b\cdot y\,:\, \bra{\lambda,y}\in E\bra{\sidx{n+1},\sidx{n+1}}}}.\end{aligned}
\end{equation}
 Let $c=\max\set{a\cdot \lambda+b\cdot y\,:\, \bra{\lambda,y}\in Q\bra{\mathcal{P},H}}$, $J^-=\set{j\in \sidx{n+1}\,:\, a\cdot\e^j+b\cdot h^{j-1} =c}=\set{j\in \sidx{n+1}\,:\, a\cdot\e^j+ b\cdot h^{j} -b\cdot c^{j-1} =c}$ and $J^+=\set{j\in \sidx{n+1}\,:\, a\cdot\e^j+b\cdot h^{j} =c}$. Then $F=\conv\bra{E\bra{J^-,J^+}}$, $b\cdot c^{j-1}=0$ for all $j\in J^+\cap J^-$, $b\cdot c^{j-1}< 0$ for all $j\in J^-\setminus J^+$ and $b\cdot c^{j-1}> 0$ for all $j\in J^+\setminus J^-$. Conversely if $F=\conv\bra{E\bra{J^-,J^+}}$ for $J^-,J^+\subseteq \sidx{n+1}$ and $b\in \Real^k$ satisfies \eqref{facelemmacond}
let $a_j=b\cdot h^j$ for $j\in J^+=\bra{ J^+\cap J^-}\cup \bra{J^+\setminus J^-}$,  $a_j=b\cdot h^{j-1}$ for $j\in J^-\setminus J^+$ and $a_j=\min\bra{b\cdot h^{j},b\cdot h^{j-1}}-1$ for $j\in \sidx{n+1}\setminus \bra{J^-\cup J^+}$. The $a$, $b$ and $F$ satisfy \eqref{facemax}.
\Halmos\endproof

\vspace{0.2in}
\proof{Proof of Proposition~\ref{sos1prophyper}.}
Throughout the proof we let $h^0=h^1$, $h^{n+1}=h^n$ as in the statement of Lemma~\ref{facelemma} so we can define $c^i=h^{i+1}-h^i$ for all $i\in \sidx{0,n}$ and it coincides with the proposition's statement for $i\in \sidx{n-1}$ and $c^0=c^{n}={\bf 0}$. We will also use the straightforward fact that $L\bra{H}=\spann\bra{\set{c^i}_{i\in \sidx{n-1}}}=\spann\bra{\set{c^i}_{j\in \sidx{0,n}}}$.

Any $\bra{\lambda,y}\in Q\bra{\mathcal{P},H}$ satisfies the equations in \eqref{sos2equalities}. Furthermore, because $L\bra{H}=\spann\bra{\set{c^i}_{i\in \sidx{n-1}}}$ we have that the dimension of the affine subspace described by \eqref{sos2equalities} is $n+\dim\bra{\set{c^i}_{i\in \sidx{n-1}}}$. Finally, by Lemma~\ref{facelemma} $\dim\bra{Q\bra{\mathcal{P},H}}=n+\dim\bra{\set{c^i}_{i\in \sidx{n-1}}}$, which shows the  statements about $\dim\bra{Q\bra{\mathcal{P},H}}$ and $\aff\bra{Q\bra{\mathcal{P},H}}$. 

Using Lemma~\ref{facelemma} we have that a face  $F\bra{J^-,J^+}=\conv\bra{E\bra{J^-,J^+}}$ for $J^-,J^+\subseteq \sidx{n+1}$ may be a facet of $Q\bra{\mathcal{P},H}$ only if $\abs{J^+\cup J^-}=n+1$ and $\dim\bra{\set{c^{j-1}}_{j\in J^+\cap J^-}}=\dim\bra{L\bra{H}}-1$ or $\abs{J^+\cup J^-}=n$ and $\dim\bra{\set{c^{j-1}}_{j\in J^+\cap J^-}}=\dim\bra{L\bra{H}}$.  

In the first option $J^+\cup J^-=\sidx{n+1}$ and $J^+\cap J^-\subsetneq \sidx{n+1}$ because $\dim\bra{\set{c^{j-1}}_{j\in \sidx{n+1}}}=\dim\bra{L\bra{H}}$. Then, $\bra{J^-\setminus J^+}\cup \bra{ J^+\setminus J^-}\neq \emptyset$ and  condition \eqref{facelemmacond} of Lemma~\ref{facelemma} holds for $b\in \Real^k\setminus \set{\bf 0}$. Furthermore, because of the first part of condition \eqref{facelemmacond} and $\dim\bra{\set{c^{j-1}}_{j\in J^+\cap J^-}}=\dim\bra{L\bra{H}}-1$ we further have that $b=s b^l$ for some $l\in \sidx{L}$ and $s\in \set{-1,1}$, and  $J^+\cap J^-=\set{j\in \sidx{n+1}\,:\, c^{j-1}\in M\bra{b^l}}$ (note  that because $c^0=c^{n}=0$ we  always have $1,n+1\in J^+\cap J^-$). If $s=1$ the second inequality in  \eqref{sos2generalfacets} for $b^l$ is satisfied at equality by all points in $E\bra{J^-,J^+}$ and strictly by all points in $E\bra{\sidx{n+1},\sidx{n+1}}\setminus E\bra{J^-,J^+}$, and hence defines $F\bra{J^-,J^+}$. Similarly, if $s=-1$ the first inequality in  \eqref{sos2generalfacets} for $b^l$ defines $F\bra{J^-,J^+}$.

In the second option we have that $\spann\bra{\set{c^{j-1}}_{j\in J^+\cap J^-}}=L(H)$ and hence the first part of condition \eqref{facelemmacond} of Lemma~\ref{facelemma} 
implies $b\in L(H)^{\perp}$. However, $\spann\bra{\set{c^{j-1}}_{j\in \sidx{n+1}}}=L(H)$ so the second part of condition \eqref{facelemmacond} implies $J^-\setminus J^+ = J^+\setminus J^-=\emptyset$ and hence  $J^-=J^+=\sidx{n+1}\setminus \set{j_0}$ for some $j_0\in J$. Noting that $c^0=c^{n}=0$ we have that condition $\spann\bra{\set{c^{j-1}}_{j\in J^+\cap J^-}}=\spann\bra{\set{c^{j-1}}_{j\in\sidx{n+1}\setminus \set{j_0}}}=L(H)$ holds if and only if $j_0\in J$ for $J$ defined in the proposition statement. Furthermore, in such case the inequality  in \eqref{sos2boundfacets} corresponding to  $j=j_0$  defines   $F\bra{J^-,J^+}$.

Hence the equations of \eqref{SOS2char} are precisely those defining $\aff\bra{Q\bra{\mathcal{P},H}}$, every facet of $Q\bra{\mathcal{P},H}$ is defined by an inequality of \eqref{SOS2char} and every inequality of \eqref{SOS2char} is facet defining for $Q\bra{\mathcal{P},H}$. Finally, the last statement follows because the two options considered for facets of $Q\bra{\mathcal{P},H}$ yield two distinct classes of facets. \Halmos\endproof

\subsection{An Embedding Formulation for SOS2 Constraints Whose Validity is Not Evident}\label{examplesection}
\begin{example} Let $\mathcal{P}:=\bra{P^i}_{i=1}^9$ be the SOS2 constraint on $\Delta^{10}$ and $H=\bigl((0, 1, 1, 1)^T,$ $(0, 1, 0, 0)^T, (0, 0, 0, 0)^T, (0, 1, 0, 1)^T, (0, 0, 0, 1)^T, (1, 0, 0, 0)^T, (1, 1, 0, 1)^T, (1, 0, 1, 1)^T, (1, 1, 1, 1)^T\bigr)$. Then   $\bra{c^i}_{i=1}^{8}=\bigl\{ (0, 0, -1, -1)^T, (0, -1, 0, 0)^T, (0, 1, 0, 1)^T, (0, -1, 0, 0)^T, (1, 0, 0, -1)^T, (0, 1, 0, 1)^T,$ $(0, -1, 1, 0)^T, (0, 1, 0, 0)^T\bigr\}$   and the set of hyperplanes spanned by them is given by  $\set{M\bra{b^l}}_{l=1}^{5}$ for $b^1=\bra{1,0,0,-1,1}^T$, $b^2=\bra{1,0,0,1}^T$, $b^3=\bra{1,-1,-1,1}^T$, $b^4=\bra{1,0,0,0}^T$ and $b^5=\bra{0,0,1,0}^T$. Finally, $\aff(H)=L(H)=\mathbb{R}^4$ and $\mathbb{R}^4=\spann\bra{\set{c^i}_{i\in \sidx{9}\setminus\set{j-1} }}$ if and only if $j\in \sidx{2,10}\setminus \set{6}$. Then \eqref{SOS2char} in this case is given by 
\begin{subequations}\label{nine}
\begin{alignat}{3}
\sum\nolimits_{j=1}^{10} \lambda_j &=1,\quad \label{nineeq}\\
\lambda_5+\lambda_6+\lambda_7+\lambda_8+\lambda_9+\lambda_{10} &\leq y_1-y_3+y_4 \\
\lambda_4+\lambda_5+\lambda_6+2\lambda_7+2\lambda_8+\lambda_9+\lambda_{10} &\geq y_1-y_3+y_4 \label{ninec}\\
\lambda_1+\lambda_5+\lambda_6+\lambda_7+2\lambda_8+2\lambda_9+2\lambda_{10} &\leq y_1+y_4 \\
\lambda_1+\lambda_2+\lambda_4+\lambda_5+\lambda_6+2\lambda_7+2\lambda_8+2\lambda_9+2\lambda_{10} &\geq y_1+y_4 \\
-\lambda_1-\lambda_2-\lambda_3+\lambda_6+\lambda_7+\lambda_8 &\leq y_1-y_2-y_3+y_4 \label{ninef}\\
-\lambda_1-\lambda_2+\lambda_5+\lambda_6+\lambda_7+\lambda_8 +\lambda_9 &\geq y_1-y_2-y_3+y_4 \\
\lambda_7+\lambda_8 +\lambda_9+\lambda_{10} &\leq y_1 \\
\lambda_6+\lambda_7+\lambda_8 +\lambda_9+\lambda_{10} &\geq y_1 \\
\lambda_1 +\lambda_9+\lambda_{10}&\leq y_3 \\
\lambda_1+\lambda_2 +\lambda_8+\lambda_9+\lambda_{10} &\geq y_3 \\
\lambda_j&\geq 0 &\quad& \forall j\in \sidx{10}\setminus\set{6},
\end{alignat}
If $y=h^4=\bra{0,1,0,1}^T$ then \eqref{nine} should enforce $\lambda_i=0$ for all $i\notin\set{4,5}$. Inequality \eqref{ninef} enforces that $\lambda_6=1$ cannot hold, but it does not force $\lambda_6=0$ as $\lambda_1=\lambda_6=1/2$ is valid for this inequality. However, this last point is infeasible for \eqref{ninec}. These two inequalities plus \eqref{nineeq} do indeed imply $\lambda_6\leq 0$ (and hence  $\lambda_6=0$ because of the lower bounds) when $y=h^4$ as adding \eqref{nineeq}, \eqref{ninec} and \eqref{ninef} yields
\[\lambda_6 \leq 1-y_2.\]
Furthermore, removing any one of these constraints allows $\lambda_6>0$ when   $y=h^4$.
\end{subequations} 
\end{example}

\subsection{Proof of Proposition~\ref{antigraylemma}}\label{omproofapendix}

To prove Proposition~\ref{antigraylemma} we need the following definition.

\begin{definition} Let $n=2^k$ for some $k\in \Int$. We say $H=\set{h^i}_{i=1}^{n}\in \mathcal{H}_{k}(n)$ is a \emph{anti-gray code}\footnote{The class of codes obtained by switching $n$ and $n-1$ in this definition is sometimes also referred to as anti-gray code.} if and only if $\sum\nolimits_{j=1}^{k} \abs{h^{2i-1}_j-h^{2i}_j}=n $ for all $i\in \sidx{ n/2}$ and $\sum\nolimits_{j=1}^{k} \abs{h^{2i}_j-h^{2i+1}_j}=n -1 $ for all $i\in \sidx{n/2-1}$.
\end{definition}

Anti-gray codes exist for all $k$ and can easily be constructed from gray codes (e.g. \cite{robinson1981counting}).  

\vspace{0.2in}
\proof{Proof of Proposition~\ref{antigraylemma}.}
Let $H$ be an anti-gray code and $c^i=h^{i+1}-h^i$ for $i\in \sidx{n-1}$. Because $H$ is an anti-gray code there exist $I\subseteq \sidx{n-1}$ with $\abs{I}=2^{k-1}$ such that $c^i\in \set{-1,1}^k$ for all $i\in I$. In addition, because $h^i\neq h^j$ for $i\neq j$ we have that $c^i\neq -c^j$ for all $i,j \in I$. Hence for all $s\in \set{-1,1}^k$ there exist $i\in I$ such that $s=c^i$ or $s=-c^i$. The result then follows from Proposition~\ref{sos1prophyper} by noting that $\set{c^i}_{i\in I}$ and $\set{\pm c^i}_{i\in I}$ span the same set of linear hyperplanes and that the number of linear hyperplanes spanned by $\set{-1,1}^k$ is equal to the number of affine hyperplanes spanned by $\set{0,1}^{k-1}$ (e.g. \cite{da2005recursivity})
\Halmos\endproof

\subsection{Proof of Proposition~\ref{UnionJackCoro}}\label{proofofUnionJackCoro}

\proof{Proof Proposition~\ref{UnionJackCoro}.} For any $H\in \mathcal{H}_k(n)$ we have $Q\bra{\mathcal{P},H}\subseteq\Delta_2^{m+1}:=\set{\lambda\in \Real^{\sidx{m+1}^2}_+\,:\, \sum_{u,v=1}^{m+1} \lambda_{\bra{u,v}}=1}$. We begin by showing that for any $H\in \mathcal{H}_k(n)$ all inequalities $\lambda_{\bra{u,v}}\geq 0$ of $\Delta_2^{m+1}$ are facet defining for $Q\bra{\mathcal{P},H}$. For that let $n=2m^2$ and $\bra{S_i}_{i=1}^{n}$ be such that $S_{2\bra{m (u-1)+(v-1)}+t}=T_{u,v}^t$ for all $u,v\in\sidx{m}$ and $t\in \set{1,2}$ so that $Q\bra{\mathcal{P},H}=\conv\bra{\bigcup_{i=1}^n P\bra{S_i}\times \set{h^i}}$. Then $\lambda_{\bra{u,v}}\geq 0$ describes a face of $Q\bra{\mathcal{P},H}$ because it is valid and is satisfied at equality for all points $\bra{\lambda,y}=\bra{\e^{\bar{u},\bar{v}},h^i}$
for $\bra{\bar{u},\bar{v}}\in S_i\setminus \set{\bra{u,v}}$. Let 
\begin{equation}\label{genericineql2222}
 a\cdot\lambda + b\cdot y\leq c 
\end{equation}
be a valid inequality of $Q\bra{\mathcal{P},H}$ that induces a facet containing the face induced by $\lambda_{\bra{u,v}}\geq 0$. Because $Q\bra{\mathcal{P},H}\subseteq\Delta_2^{m+1}\cap \aff\bra{H}$, without loss of generality we may assume $c=0$ and $b\in L\bra{H}$ by possibly adding multiples of $\sum_{u,v=1}^{m+1} \lambda_{\bra{u,v}}=1$ and the equations defining $\aff\bra{H}$. For any $i,j\in \sidx{n}$, $\bra{\underline{u},\underline{v}}\in S_i \setminus \set{\bra{u,v}}$ and $\bra{\overline{u},\overline{v}}\in S_j \setminus \set{\bra{u,v}}$ there exist $\bra{i_l}_{l=1}^r\subseteq \sidx{n}$  and $\bra{\bra{u_l,v_l}}_{l=1}^{r-1} \subseteq \sidx{n+1}^2\setminus \set{\bra{u,v}}$ such that $i_1=i$, $i_r=j$, $\bra{\underline{u},\underline{v}}=\bra{u_{i_1},v_{i_1}}$, $\bra{\overline{u},\overline{v}}=\bra{u_{i_r},v_{i_r}}$ and $\bra{u_l,v_l}\in S_{i_l}\cap S_{i_{l+1}}$ for all $l\in \sidx{r-1}$. Because $\bra{\lambda,y}=\bra{\e^{\bra{u_l,v_l}},h^l}$ and $\bra{\lambda,y}=\bra{\e^{\bra{u_l,v_l}},h^{l+1}}$ satisfy \eqref{genericineql2222} at equality for all $l\in \sidx{r-1}$, we have $b\cdot h^i=b\cdot h^j$ and $a_{\bra{\underline{u},\underline{v}}}=a_{\bra{\overline{u},\overline{v}}}=b\cdot h^i$. Because the first identity holds for all $i,j\in \sidx{n}$ we have  $b\cdot\bra{h^i-h^1}=0$ for all $i\in\sidx{n}$, which together with $b\in L\bra{H}$ implies $b=\bf 0$. Similarly, combining   $b=\bf 0$ and the fact that the second identity holds for all $\bra{\underline{u},\underline{v}},\bra{\overline{u},\overline{v}}\neq \bra{u,v}$ 
we obtain $ a_{\bra{\bar{u},\bar{v}}}=0$ for all $\bra{\bar{u},\bar{v}}\neq  \bra{u,v}$. Finally,  validity of \eqref{genericineql2222} implies $a_{\bra{u,v}}\leq 0$ and being facet defining further implies $a_{\bra{u,v}}<0$. Then \eqref{genericineql2222} is a positive multiple of $\lambda_{\bra{u,v}}\geq 0$  and hence $\lambda_{\bra{u,v}}\geq 0$ is facet defining.   

The result on $\lambda_{\bra{u,v}}\geq 0$ shows that $\bra{\sqrt{n/2}+1}^2\leq \size_B\bra{Q\bra{\mathcal{P},H}} \leq \mmc\bra{\mathcal{P}} $  for any $H\in \mathcal{H}_k(n)$.

The sizes for the unary encoded formulation comes  from Proposition~10 in \cite{lee01}, the comments before its statement and Lemma~\ref{formulationLemma}. The existence and sizes for the binary encoded formulation come from  the proof of Theorem~1 in \cite{Modeling-Disjunctive-Constraints-FULL}, Lemma~\ref{formulationLemma}, the lower bound on $\size_B\bra{Q\bra{\mathcal{P},H}}$,  noting that the bounds on the binary variables of the formulation from \cite{Modeling-Disjunctive-Constraints-FULL} are redundant (cf. the logarithmic formulation for SOS2) and that none of the non-bound inequalities of this formulation are redundant.\Halmos\endproof

\ACKNOWLEDGMENT{The author would like to thank the review team including an anonymous associate editor and two anonymous referees for their thoughtful and constructive comments, which significantly improved the exposition of the paper. This research was partially supported by the National Science Foundation under grant CMMI-1351619.}

\bibliographystyle{ormsv080}
%\bibliography{references}

\end{document}